\documentclass[10pt]{amsart}

\usepackage[encapsulated]{CJK}
\usepackage[utf8]{inputenc}
\usepackage{lmodern}
\usepackage{graphicx}   
\usepackage{bm}
\usepackage[english]{babel}
\usepackage{amsmath,amssymb,amsthm}
\usepackage{mathcomp}
\usepackage{setspace}
\usepackage{hyperref}
\usepackage{xcolor}
\usepackage[a4paper]{geometry}

\usepackage{booktabs}
\let\rho\varrho
\theoremstyle{remark}
\newtheorem{remark}{Remark}
\newcommand{\dxdy}{\,\mathrm{d}x\mathrm{d}y}
\newcommand{\dxideta}{\,\mathrm{d}\xi\mathrm{d}\eta}
\newcommand{\dxi}{\,\mathrm{d}\xi}
\newcommand{\dtau}{\,\mathrm{d}\tau}
\newcommand{\dtaudsigma}{\,\mathrm{d}\tau\mathrm{d}\varsigma}
\newcommand{\dx}{\,\mathrm{d}x}
\numberwithin{equation}{section}
\definecolor{viola}{rgb}{0.4,0.,0.9}
\newcommand{\revCom}[1]{\textcolor{black}{#1}}
\begin{document}
%
\title[Multi-element SIAC filter for shock capturing DGSEM]{Multi-element SIAC filter for shock capturing applied to high-order discontinuous Galerkin spectral element methods}

\author[bohm]{Marvin Bohm$^{1}$}
\author[schermeng]{Sven Schermeng$^{1}$}
\address{$^{1}$Mathematical Institute, University of Cologne, Weyertal 86-90, 50931 Cologne, Germany}
\author[winters]{Andrew R. Winters$^{2,*}$}
\address{$^{2}$Department of Mathematics; Computational Mathematics, Link\"{o}ping University, SE-581 83 Link\"{o}ping, Sweden}
\author[gassner]{Gregor J. Gassner$^{3}$}
\address{$^{3}$Department for Mathematics and Computer Science; Center for Data and Simulation Science, University of Cologne, Weyertal 86-90, 50931 Cologne, Germany}
\author[jacobs]{Gustaaf B. Jacobs$^{4}$}
\address{$^4$Department of Aerospace Engineering, San Diego State University, 5500 Campanile Drive, San Diego, CA 92182, United States}

\maketitle

\begin{abstract}
We build a multi-element variant of the smoothness increasing accuracy conserving (SIAC) shock capturing technique proposed for single element spectral methods by Wissink et al. \cite{wissink2018shock}. In particular, the baseline scheme of our method is the nodal discontinuous Galerkin spectral element method (DGSEM) for approximating the solution of systems of conservation laws. It is well known that high-order methods generate spurious oscillations near discontinuities which can develop in the solution for nonlinear problems, even when the initial data is smooth. We propose a novel multi-element SIAC filtering technique applied to the DGSEM as a shock capturing method. We design the SIAC filtering such that the numerical scheme remains high-order accurate and that the shock capturing is applied adaptively throughout the domain. The shock capturing method is derived for general systems of conservation laws. We apply the novel SIAC filter to the two-dimensional Euler and ideal magnetohydrodynamics (MHD) equations to several standard test problems with a variety of boundary conditions. 
\end{abstract}
\vspace{0.5cm}

\section{Introduction}
\revCom{Systems of nonlinear hyperbolic conservation laws cover a wide range of physical flow problems, e.g. modeled by the Euler or magneto-hydrodynamic (MHD) equations. Such flow phenomena are particularly interesting as they describe diverse physical processes like gas dynamics in chemical processes or plasma interactions in space, respectively \cite{kundu2008,Landau1959}. It is well-known that the solution of nonlinear hyperbolic systems can develop discontinuities, e.g. shocks, in finite time, regardless of the smoothness of the initial data \cite{evans2010partial}. Due to the complexity of nonlinear systems, we rely on numerical methods to approximate the solutions to such problems. }

\revCom{For low-order spatial approximations, like finite volume methods, discontinuous solutions are unproblematic because their inherent amount of numerical dissipation regularizes discontinuities naturally. However, for high-order numerical methods spurious oscillations near discontinuities, i.e. {Gibbs phenomenon} \cite{gottlieb1997gibbs}, arise. These unphysical overshoots might lead to unphysical solution states, e.g. negative density or pressure. Over the decades, many counter mechanisms have been developed to control overshoots and stabilize high-order approximations in shocked regions. Altogether, these methods can be subdivided into three main categories, which are slope limiters \cite{arora1997well,balsara2000monotonicity,garnier1999use,harten1983high,harten1987uniformly,shu1998essentially}, artificial dissipation techniques \cite{chaudhuri2017,hartmann2006adaptive,hughes1986new,jameson1981numerical,persson2006sub}, and solution filters \cite{alvarez1994signal,chaudhuri2017,mirzaee2011smoothness,ryan2015one}.}

\revCom{In this work, we consider a shock capturing method that uses the global smoothness increasing accuracy conserving (SIAC) filter, a filter that has received much attention in the context of postprocessing data produced by discontinuous Galerkin approximations, see e.g. \cite{mirzaee2011smoothness,ryan2015one,van2011position}. The SIAC filter increases smoothness by convolving the approximate solution with an appropriate smooth kernel function, e.g. B-Splines \cite{ji2015smoothness,ryan2015one} or Dirac-delta polynomial approximations, e.g. \cite{suarez2017regularization}. The accuracy conservation is more technical and is related to the fundamental building block of the discontinuous Galerkin solution space(s) and its solution ansatz \cite{cockburn2003}. The filter kernel is designed to reproduce certain polynomials orders by convolution, for example $m$. Consider a DG approximation that uses a polynomial solution space of $N$. If the filter is designed to recover a larger family of polynomial orders, i.e. $m > N$, then the filter conserves the accuracy of the approximation. If the filter recovers a smaller family of polynomial orders, $m < N$, then the accuracy of the overall approximation is bound by the filter accuracy. Typically, such SIAC filters were designed to obtain super-convergence in a post-processing step by a convolution of the numerical approximation against a specifically designed kernel function once at the final time \cite{bramble1977,cockburn2003,ji2015smoothness,li2019smoothness,ryan2015one,steffen2008}. However, recent work has applied the SIAC filter as a shock capturing and/or regularization of general discontinuities strategy during the computation of the approximate solution for global spectral collocation methods \cite{suarez2017regularization,wissink2018shock}. For such spectral methods, the SIAC filter is suitable to apply in shocked regions because said filter can recover full accuracy away from shocks, see e.g \cite{ji2015smoothness}. The filter for global spectral collocation methods is constructed with a Dirac-delta kernel sequence determined by two conditions that control the number of vanishing moments and the smoothness. With the discrete SIAC filter at hand, the shock regularization of a global approximation is then performed by a convolution of the solution with the high-order Dirac-delta kernels in every time step. In practice, the filtering procedure reduces to a simple matrix-vector multiplication and, thus, allows for an efficient and simple implementation.}

\revCom{The main contribution of this work is to extend the global filtering technique to element-wise approximations within a nodal discontinuous Galerkin (DG) method. We exploit the weak coupling of the discontinuous Galerkin spectral element method (DGSEM) at element interfaces to design a \textit{multi-element filter}.} Consequently, we convolve the polynomial approximations of one element and its nearest neighbor's solutions with Dirac-delta kernels instead of the global representation of the solution. As we will point out in the derivation of this {multi-element filter}, it also recasts to \textit{locally} applied matrix-vector formulations. We present the filtering matrices for one-dimensional and two-dimensional Cartesian DG discretizations. Moreover, we construct the multi-element SIAC filter such that the numerical scheme remains high-order accurate in smooth regions and that the shock capturing is applied adaptively throughout the domain. The latter we obtain by implementing a shock indicator, which is defined by the filtered solution itself. According to this indicator, we replace the DG approximation by the filtered one in oscillatory regions and, additionally, even introduce a smooth transition area, in which we use a convex combination of the filtered and unfiltered solutions.

The outline of this work is as follows: We begin in Sec.~\ref{Sec:DG} with a construction of the DG method on two-dimensional Cartesian meshes. Next, in Sec.~\ref{sec:SIAC} we provide the SIAC filtering routines, where we first discuss the global filter, before we extend it to the multi-element variant as well as two spatial dimensions. Furthermore, we broach the issue of adaptive filtering and conservation properties in the same section. Finally, we provide several numerical benchmark tests in order to verify the applicability of the novel filter to shock problems for the two-dimensional Euler and ideal MHD equations in Sec.~\ref{Sec:NumTest}. Lastly, Sec.~\ref{Sec:Conc} gives concluding remarks and an outlook on possible further research projects.

\section{Discontinuous Galerkin spectral element method}\label{Sec:DG}

Throughout this work we consider the solution of hyperbolic conservation laws in two spatial dimensions which take the general form
\begin{equation}\label{eq:consLawGeneral}
u_t + f(u)_x + g(u)_y = 0,
\end{equation}
in a square domain $\Omega = \left[x_L,x_R\right] \times \left[y_B,y_T\right] \subset \mathbb{R}^2$ with appropriate boundary conditions and an initial solution $u(x,y,0) = u_0(x,y)$. Here $u$ is a conserved variable and $f,g$ are the nonlinear fluxes. We take \eqref{eq:consLawGeneral} as the prototype equation for the conserved solution variables such as density or momentum. This simplifies the discussion and derivations for the SIAC filtering in Sec. \ref{sec:SIAC}. The discussion extends naturally to a system of hyperbolic conservations laws such as the Euler equations.

We first provide an overview for the derivation of the nodal DGSEM on Cartesian meshes. Complete details can be found in the book by Kopriva \cite{Kopriva:2009nx}. We derive the DGSEM from the weak form of the conservation law \eqref{eq:consLawGeneral}. As such, we multiply by an arbitrary discontinuous $L_2(\Omega)$ test function $\varphi$ and integrate over the domain
\begin{equation}\label{eq:weakForm}
\int\limits_{\Omega}\left({u}_t+{f}_x+{g}_y\right)\varphi\dxdy = {0},
\end{equation}
where we suppress the ${u}$ dependence of the nonlinear fluxes. We subdivide the domain $\Omega$ into $N_Q$ non-overlapping quadratic elements
\begin{equation}
Q_n = \left[x_{n,1},x_{n,2}\right]\times[y_{n,1},y_{n,2}],\quad n = 1,\ldots,N_Q.
\end{equation}
For the present discussion we make the simplifying assumption that all elements have the same size, i.e. $\Delta x = x_{n,2}-x_{n,1}$and $\Delta y = y_{n,2}-y_{n,1}$ for all $n=1,\ldots,N_Q$. This divides the integral over the whole domain into the sum of the integrals over the elements. So, each element contributes
\begin{equation}\label{eq:weakForm2}
\int\limits_{Q_n}\left({u}_t+{f}_x+{g}_y\right)\varphi\dxdy = {0},\quad n = 1,\ldots N_Q,
\end{equation}
to the total integral. Next, we create a transformation between the reference element $Q_0=[-1,1]^2$ and each element, $Q_n$. For rectangular meshes we create mappings $(X_n,Y_n):Q_0 \rightarrow Q_n$ such that $\left(X_n(\xi),Y_n(\eta)\right) = (x,y)$ are
\begin{equation}
 X_n(\xi) = x_{n,1} + \frac{\xi +1}{2}\Delta x,\quad  Y_n(\eta) = y_{n,1} + \frac{\eta +1}{2}\Delta y,
\label{eq:mapping}
\end{equation}
for $n=1,\ldots,N_Q$. Under the transformation \eqref{eq:mapping} the conservation law in physical coordinates \eqref{eq:consLawGeneral} becomes a conservation law in reference coordinates \cite{Kopriva:2009nx}
\begin{equation}
J{u}_t + \tilde{{f}}_{\xi} + \tilde{{g}}_{\eta} = {0},
\end{equation}
where
\begin{equation}
J=\frac{\Delta x\Delta y}{4},\quad\tilde{{f}} = \frac{\Delta y}{2}{f},\quad\tilde{{g}} = \frac{\Delta x}{2}{g}.
\end{equation}

We select the test function $\varphi$ to be a piecewise polynomial of degree $N$ in each spatial direction
\begin{equation}\label{eq:testFunction}
\varphi^{Q_n} = \sum_{i=0}^N\sum_{j=0}^N\varphi_{ij}^{Q_n}\psi_i(\xi)\psi_j(\eta),
\end{equation}
on each spectral element $Q_n$, but do not enforce continuity at the element boundaries. The interpolating Lagrange basis functions are defined by
\begin{equation}
\psi_i(\xi) = \prod_{\stackrel{j = 0}{j \neq i}}^N \frac{\xi-\xi_j}{\xi_i-\xi_j} \quad\text{for}\quad i=0,\ldots,N,
\label{eq:Lagrange}
\end{equation}
with a similar definition in the $\eta$ direction. The values of $\varphi^{Q_n}_{ij}$ on each element $Q_n$ are arbitrary and linearly independent, therefore the formulation \eqref{eq:weakForm2} is
\begin{equation}
\int\limits_{Q_0}\left(J{u}_t+\tilde{{f}}_{\xi}+\tilde{{g}}_{\eta}\right)\psi_i(\xi)\psi_j(\eta)\dxideta = {0},
\end{equation}
where $i,j=0,\ldots,N$.

For the DGSEM we approximate the conservative variable $u$ and the contravariant fluxes $\tilde f, \tilde g$ with polynomial interpolants of degree $N$ in each spatial direction written in Lagrange form on each element $Q_n$, e.g.,
\begin{equation}
\begin{aligned}
u(x,y,t)|_{Q_n} &= u^{Q_n}(\xi,\eta,t) \approx \sum_{i=0}^N\sum_{j=0}^N u^{Q_n}_{ij}(t) \psi_i(\xi) \psi_j(\eta),\\[0.05cm]
\tilde f(u(x,y,t))|_{Q_n} &= \tilde{f}^{Q_n}(\xi,\eta,t) \approx \sum_{i=0}^N\sum_{j=0}^N \tilde{f}^{Q_n}_{ij}(t) \psi_i(\xi) \psi_j(\eta).
\end{aligned}
\label{DG-approx}
\end{equation} 
This implies that the global representation of the solution $u$ is the union of these piecewise polynomials
\begin{equation}
u(x,y,t) \approx \bigcup_{n = 1}^{N_Q} u^{Q_n}(\xi,\eta,t).
\label{DG-global}
\end{equation}

Next, we use integration-by-parts to move derivatives off the nonlinear fluxes and onto the test function, which generates surface and volume contributions. We resolve the discontinuities between elements at the surface by introducing Lax-Friedrichs numerical flux functions $f^*$ and $g^*$. We apply integration-by-parts a second time to move derivatives back onto the fluxes. For the nodal DGSEM any integrals present in the variational formulation are approximated with $N+1$ Legendre-Gauss-Lobatto (LGL) quadrature nodes and weights, e.g.,
\begin{equation}\label{eq:quadrature}
\begin{split}
\int\limits_{Q_0} Ju^{Q_n}_t\psi_i(\xi)\psi_j(\eta)\dxideta & \approx \!\!\!\sum_{n,m=0}^N\left(\sum_{p,q=0}^NJ\left(u^{Q_n}_t\right)_{pq}\psi_p(\xi_n)\psi_q(\eta_m)\right)\psi_i(\xi_n)\psi_j(\eta_m)\omega_n \omega_m
\\ & = J\left(u^{Q_n}_t\right)_{ij}\omega_i \omega_j,
\end{split}
\end{equation}
for each element $n=1,\ldots,N_Q$ and where $\left\{\xi_i\right\}_{i=0}^N, \left\{\eta_j\right\}_{j=0}^N$ are the LGL quadrature nodes and $\left\{\omega_i\right\}_{i=0}^N,\left\{\omega_j\right\}_{j=0}^N$ are the LGL quadrature weights. Further, we \textit{collocate} the interpolation and quadrature nodes which enables us to exploit that the Lagrange basis functions \eqref{eq:Lagrange} are discretely orthogonal and satisfy the Kronecker delta property, i.e., $\psi_j(\xi_i) = \delta_{ij}$ with $\delta_{ij} = 1$ for $i=j$ and $\delta_{ij}=0$ for $i\neq j$ to simplify \eqref{eq:quadrature}. Due to the polynomial approximation \eqref{DG-approx} any derivatives fall on the Lagrange basis functions. These are approximated at high-order with the standard differentiation matrix \cite{Kopriva:2009nx}
\begin{equation}
D_{ij} = \psi_j'(\xi)\bigg|_{\xi=\xi_i},\quad i,j = 0,\ldots,N.
\end{equation}

From these steps, we build the standard semi-discrete formulation of the strong-form DGSEM. We write the scheme in index notation as
\begin{equation}\label{2D-DG-semi}
\begin{aligned}
J\left(u^{Q_n}_t\right)_{ij} = -&\left(\frac{\delta_{iN}}{\omega_N}\left[\tilde f^{*}(1,\eta_j) - \tilde f^{Q_n}_{ij}\right] - \frac{\delta_{i0}}{\omega_0}\left[\tilde f^{*}(-1,\eta_j) - \tilde f^{Q_n}_{ij}\right]+\sum_{m=0}^N D_{im}\tilde f^{Q_n}_{mj}\right.\\
&\;\;\;\left.\frac{\delta_{Nj}}{\omega_N}\left[\tilde g^{*}(\xi_i,1) - \tilde g^{Q_n}_{ij}\right] - \frac{\delta_{0j}}{\omega_0}\left[\tilde g^{*}(\xi_i,-1) - \tilde g^{Q_n}_{ij}\right]+\sum_{m=0}^N D_{jm}\tilde g^{Q_n}_{im}\right),
\end{aligned}
\end{equation}
where $i,j=0,\ldots,N$ and $n=1,\ldots,N_Q$.

The semi-discrete formulation \eqref{2D-DG-semi} on each element is integrated in time with an explicit five stage, fourth order Runge-Kutta method of Carpenter and Kennedy \cite{Kennedy1994}. We select a stable explicit time step with an appropriate CFL condition which is equation and resolution dependent.

\section{SIAC Filter}\label{sec:SIAC}

In this section we present the SIAC filter for a single domain and then discuss its extension and implementation into a multi-element DGSEM framework.


\subsection{Single element filter}

\begin{figure}[htbp!]
\begin{center}
\includegraphics[width=0.8\textwidth]{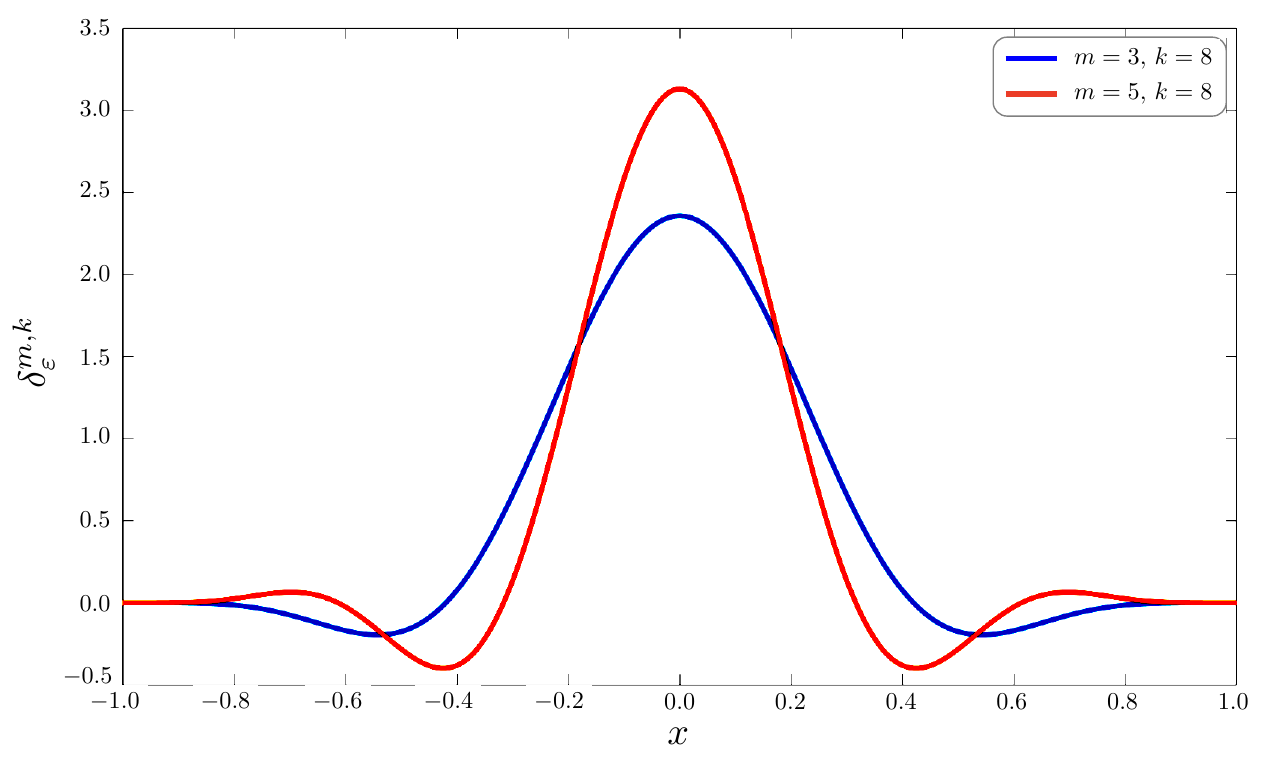}
\caption{Visualization of the Dirac-delta kernel for $\varepsilon=1$.}
\label{fig:deltakernel}
\end{center}
\end{figure}

In \cite{wissink2018shock} a SIAC filter was developed for a single domain spectral method. To define the global method we begin by introducing the delta sequence
\begin{equation}
\delta^{m,k}_\varepsilon(x) = \begin{cases} \frac{1}{\varepsilon} P^{m,k}\left(\frac{x}{\varepsilon}\right) ~~ & \left|x\right|\leq\varepsilon \\ 0 & \left|x\right|>\varepsilon \end{cases},
\label{delta}
\end{equation}
that is built from the polynomial $P^{m,k}(x)$, which is uniquely determined by the following conditions:
\begin{align}
\int\limits_{-1}^1 P^{m,k}(\xi) \dxi & = 1 \label{moment0}, \\ 
\left(P^{m,k}\right)^{(i)}(\pm 1) & = 0\quad \text{for}\quad i=0,\ldots,k, \\
\int\limits_{-1}^1 \xi^i P^{m,k}(\xi) \dxi & = 0\quad \text{for}\quad i=1,\ldots,m. \label{moments}
\end{align} 
In Fig.~\ref{fig:deltakernel} we illustrate the polynomial approximation of the delta kernel \eqref{delta} with $(m,k) = (3,8)$, $(m,k) = (5,8)$, and scaling parameter $\varepsilon=1$.

According to the SIAC filtering strategy \cite{ryan2015one,suarez2017regularization,wissink2018shock} we regularize the solution produced by the numerical scheme with the manipulation
\begin{equation}\label{start}
\resizebox{\textwidth}{!}{$
 \tilde{u}(x,t)  = \int\limits_{x-\varepsilon}^{x+\varepsilon} u(\tau,t) \delta^{m,k}_{\varepsilon} (x-\tau) \dtau
            \approx \int\limits_{x-\varepsilon}^{x+\varepsilon} \left[\sum_{i=0}^N u_i(t) \psi_i(\tau)\right] \delta^{m,k}_{\varepsilon} (x-\tau) \dtau
					  = \sum_{i=0}^N u_i(t) \int\limits_{x-\varepsilon}^{x+\varepsilon}  \psi_i(\tau) \delta^{m,k}_{\varepsilon} (x-\tau) \dtau.
$}
\end{equation}
For compact notation we introduce the filter matrix $\Phi$ and approximate its values with LGL quadrature by mapping the corresponding integration area $\left[x-\varepsilon,x+\varepsilon\right]$ into the reference element $E=\left[-1,1\right]$ 
\begin{equation}
\label{FilterMatrixSingle}
\Phi_{ij}  = \int\limits_{x_i-\varepsilon}^{x_i+\varepsilon} \psi_j(\tau) \delta^{m,k}_{\varepsilon} (x_i-\tau) \dtau
				   = \varepsilon \int\limits_{-1}^1 \psi_j(\varepsilon x+x_i) \delta^{m,k}_{\varepsilon} (\varepsilon x) \dx
					 \approx \varepsilon \sum_{\nu=0}^{N^\ast} \omega_\nu \psi_j(\varepsilon x_\nu + x_i) \delta^{m,k}_{\varepsilon} (\varepsilon x_\nu),
\end{equation}
where $\left\{x_\nu\right\}_{\nu=0}^{N^\ast}$ and $\left\{\omega_\nu\right\}_{\nu=0}^{N^\ast}$ are the LGL quadrature points and weights for $N^\ast = 2\left(\frac{m}{2}+k+1\right)$ to maintain the desired high-order accuracy of the approximation  \cite{suarez2017regularization}. Moreover, we choose 
\begin{equation}
\varepsilon=\cos\left(\frac{\pi(\frac{N-N_d}{2})}{N}\right),
\label{support}
\end{equation}
with $N_d$ determined empirically to ensure stable conserving results \cite{suarez2017regularization,wissink2018shock}. \revCom{The selection of the scaling parameter $\varepsilon$ is related to the quadrature accuracy of the integral \eqref{FilterMatrixSingle}. As was shown in \cite{suarez2017regularization}, an arbitrary choice of $\varepsilon$ can lead to sub-optimal accuracy of the quadrature. However, depending on the number of LGL nodes used for the DG approximation the support of the kernel function \eqref{support} must be adjusted through the parameter $N_d$. In practice, for a fixed number of LGL nodes $N+1$, selecting different values of $\varepsilon$ changes the accuracy of the filter as well as the solution quality because of possible excessive smearing effects, see \cite{suarez2017regularization,wissink2018shock} for details.}

We can express the filtering process in terms of a matrix vector multiplication, i.e.
\begin{equation}
\tilde{\mathbf{u}}(t) = \mathbf{\Phi} \mathbf{u}(t).
\end{equation}
For the global SIAC filtering technique we must address how the filter matrix is applied at the physical boundaries of the domain. However, at the physical boundaries no $\varepsilon$-stencils are defined. Thus, oscillations caused by shocks as well as by re-interpolation (Runge phenomena) cannot be smoothed in these areas. In the original approach for the global collocation the affected parts of the discretization are set to the analytical solution \cite{wissink2018shock}. Using a local version of the filter we can avoid identifying interior points by an analytical reference solution, as discussed in the next section.

We also note that, by construction, the filter conserves mass solely for polynomial data of degree up to $m$, which especially in a global collocation method is difficult to realize. Thus, small conservation errors might be introduced by applying the filter matrix \eqref{FilterMatrixSingle}.

\subsection{Multi-element filter}

For the case with multiple Cartesian elements and DGSEM, we begin, again, with the one-dimensional case and then apply the tensor product decoupling of the DGSEM to move to higher spatial dimensions. In contrast to the previous section we now have DG solutions defined on multiple elements in a mesh ($N_Q>1$). However, we still want to apply the smoothing matrix $\Phi$ locally to the solution element-wise. It is, hence, necessary to determine how to couple the filtering across element interfaces to determine a \textit{multi-element} SIAC technique.

To do so, we begin with \eqref{start}, where, in one spatial dimension, we know that the solution $u$ is a union of piecewise polynomials over all elements $N_Q$
\begin{equation}\label{eq:multielement1}
\begin{aligned}
 \tilde{u}(x,t)  = \int\limits_{x-\varepsilon}^{x+\varepsilon} u(\tau,t) \delta^{m,k}_{\varepsilon} (x-\tau) \dtau &= \int\limits_{x-\varepsilon}^{x+\varepsilon} \left[\sum_{\ell=1}^{N_Q} u^{Q_\ell}(\tau,t)\right] \delta^{m,k}_{\varepsilon} (x-\tau) \dtau \\
					 & = \sum_{\ell=1}^{N_Q} \int\limits_{x-\varepsilon}^{x+\varepsilon} u^{Q_\ell}(\tau,t) \delta^{m,k}_{\varepsilon} (x-\tau) \dtau.
\end{aligned}
\end{equation} 
Next, we focus on one physical node $x_i := X_n(\xi_i)$ within one element $Q_n$ and define the following sets
\begin{align}
Q^\varepsilon_{i,n} & := \left[x_i-\varepsilon,x_i+\varepsilon\right] \cap Q_n, \\
Q^\varepsilon_{i,n-1} & := \left[x_i-\varepsilon,x_i+\varepsilon\right] \cap Q_{n-1}, \\
Q^\varepsilon_{i,n+1} & := \left[x_i-\varepsilon,x_i+\varepsilon\right] \cap Q_{n+1}.
\end{align}
Since $\varepsilon$ is sufficiently small, we assume that the $\varepsilon$-stencil is imbedded in these three sets and thus
\begin{equation}
\begin{aligned}
 \tilde{u}_i^{Q_n}(t) & = \sum_{\ell=1}^{N_Q} \int\limits_{x_i-\varepsilon}^{x_i+\varepsilon} u^{Q_\ell}(\tau,t) \delta^{m,k}_{\varepsilon} (x_i-\tau) \dtau \\
  & = \int\limits_{Q^\varepsilon_{i,n}} u^{Q_n}(\tau,t) \delta^{m,k}_{\varepsilon} (x_i-\tau) \dtau + \int\limits_{Q^\varepsilon_{i,n-1}} u^{Q_{n-1}}(\tau,t) \delta^{m,k}_{\varepsilon} (x_i-\tau) \dtau \\
	& \qquad+ \int\limits_{Q^\varepsilon_{i,n+1}} u^{Q_{n+1}}(\tau,t) \delta^{m,k}_{\varepsilon} (x_i-\tau) \dtau.
\end{aligned}
\end{equation}
If we now define similar sets for the corresponding LGL node
\begin{align}
E^\varepsilon_i & := \left[\xi_i-\varepsilon,\xi_i+\varepsilon\right] \cap \left[-1,1\right], \\
E^\varepsilon_{i,l} & := \left[\xi_i-\varepsilon,\xi_i+\varepsilon\right] \cap \left[\xi_i-\varepsilon,-1\right], \\
E^\varepsilon_{i,r} & := \left[\xi_i-\varepsilon,\xi_i+\varepsilon\right] \cap \left[1,\xi_i+\varepsilon\right], 
\end{align}
we can transform everything to reference space again and obtain
\begin{equation}\label{mappedFilter}
\begin{aligned}
 \tilde{u}_i^{Q_n}(t) & = \varepsilon \int\limits_{E^\varepsilon_i} \sum_{j=0}^N \psi_j(\varepsilon x+\xi_i) u_j^{Q_n}(t) \delta^{m,k}_{\varepsilon} (\varepsilon x) \dx + \varepsilon \int\limits_{E^\varepsilon_{i,l}} \sum_{j=0}^N \psi_j(\varepsilon x+\xi_i-2) u_j^{Q_{n-1}}(t) \delta^{m,k}_{\varepsilon} (\varepsilon x) \dx \\
& \qquad+ \varepsilon \int\limits_{E^\varepsilon_{i,r}} \sum_{j=0}^N \psi_j(\varepsilon x+\xi_i+2) u_j^{Q_{n+1}}(t) \delta^{m,k}_{\varepsilon} (\varepsilon x) \dx.
\end{aligned}
\end{equation}
Note, that we shift the arguments of the lagrange basis functions in the left and right elements by $\pm2$ to guarantee the correct evaluation points. We can write \eqref{mappedFilter} in compact notation by applying a modified $(N\!+\! 1)\times 3(N\!+\! 1)$ smoothing matrix to the solution, i.e.
\begin{equation}
   \tilde{\mathbf{u}}^{Q_n}(t) = \underbrace{\begin{pmatrix} \mathbf{\Phi}^{n-1} & \mathbf{\Phi}^n & \mathbf{\Phi}^{n+1} \end{pmatrix}}_{=\mathbf{\Phi}^{Q_n}_{loc}} \begin{pmatrix} \mathbf{u}^{Q_{n-1}}(t) \\ \mathbf{u}^{Q_n}(t) \\ \mathbf{u}^{Q_{n+1}}(t) \end{pmatrix},
\end{equation}
with the matrices
\begin{align}
 \Phi^n_{i,j} & = \varepsilon \int\limits_{E^\varepsilon_i} \psi_j(\varepsilon x+\xi_i) \delta^{m,k}_{\varepsilon} (\varepsilon x) \dx, \label{eq:multiFilterInt1}\\
 \Phi^{n-1}_{i,j} & = \varepsilon \int\limits_{E^\varepsilon_{i,l}} \psi_j(\varepsilon x+\xi_i-2) \delta^{m,k}_{\varepsilon} (\varepsilon x) \dx, \\
\Phi^{n+1}_{i,j} & =  \varepsilon \int\limits_{E^\varepsilon_{i,r}} \psi_j(\varepsilon x+\xi_i+2) \delta^{m,k}_{\varepsilon} (\varepsilon x) \dx.\label{eq:multiFilterInt3}
\end{align}
Again, we evaluate these integrals by a Legendre-Gauss-Lobatto quadrature with $N^\ast = 2\left(\frac{m}{2}+k+1\right)$ points as in \eqref{FilterMatrixSingle}. 
\revCom{Because the two dimensional, multi-element filter is created through a tensor product of the one-dimensional filters, we select the same $\varepsilon$ \eqref{support} in each of the integrals \eqref{eq:multiFilterInt1}-\eqref{eq:multiFilterInt3}.} 

Note, that since the neighboring elements only enter in the smoothing matrix at grid points near to the element boundaries, $\mathbf{\Phi}^{n-1}$ and $\mathbf{\Phi}^{n+1}$ are block matrices with mostly zero entries, especially when $N$ is large. In particular, only the first several rows of $\mathbf{\Phi}^{n-1}$ and the last several rows of $\mathbf{\Phi}^{n+1}$ are non-zero. We see that the multi-element SIAC filtering process is not entirely local to element $Q_n$; however, we only need solution information from its direct neighbors in the mesh. 

An additional advantage in the design of this multi-element filtering technique is the treatment of element located at physical boundaries. We already noted that there is no $\varepsilon$-stencil defined in these boundary areas. Thus, we cannot apply the filter. However, from the multi-element technique we can introduce \textit{ghost elements}, in which we can define a consistent solution depending on the physical boundary condition, e.g., to reflecting wall or Dirichlet values. This procedure removes Runge phenomena from the solution without the need to identify interior points by analytical values \cite{wissink2018shock}.

We can use the locally filtered solution as a \textit{shock detector} to adaptively apply the multi-element filter only in elements where it is necessary. To do so, we define an indicator to measure the difference between the filtered and unfiltered solutions
\begin{equation}
  \epsilon_n := \underset{i=0,\ldots,N}{\max}\left|\mathbf{u}^{Q_n}_i - \tilde{\mathbf{u}}^{Q_n}_i \right|.
\label{indic}
\end{equation}
Next, we normalize this indicator with respect to the polynomial order and the number of elements and check in each element $n=1,\ldots,N_Q$, if
\begin{equation}
  \frac{\epsilon_n}{(N+1)N_Q} > \texttt{TOL},
	\label{adaptive}
\end{equation}
for a given user defined tolerance \texttt{TOL}. If this condition is fulfilled, we replace the current element solution with the filtered solution. Otherwise the approximation is deemed to be sufficiently smooth and no filtering is applied. In order to extend this adaptation for systems of conservation laws we can use single variables to compute $\epsilon_n$, e.g., the density or pressure for the Euler equations.
\revCom{\begin{remark}[Adaptive filtering]
The filtered solution acts as a self-contained shock detector because of the constraints used to construct the SIAC filter \eqref{moment0}-\eqref{moments}. In particular, the filter is designed to recover polynomial orders up to $m$. Therefore, in smooth regions of the flow the approximate solution and the filtered solution will be nearly the same, i.e., the indicator error \eqref{adaptive} will be small. However, near discontinuities the approximate solution will contain large, spurious overshoots whereas the overshoots of the filtered solution will be considerably reduced, making \eqref{adaptive} large.
\end{remark}}

Moreover, for convenience, we define
\begin{equation}
\sigma_n = \log_{10} \left(\epsilon_n\right),
\end{equation}
and introduce a transition area between two tolerance levels, $\sigma_{\min}\leq\sigma_n\leq \sigma_{\max}$, to smoothly blend the filtered and unfiltered solutions. As such we define a parameter $0\leq\lambda\leq 1$ and then define the updated solution on a given element $Q_n$ to be a convex combination of the two solutions
\begin{equation}
\hat{\mathbf{u}}^{Q_n} = \lambda \tilde{\mathbf{u}}^{Q_n} + (1-\lambda) \mathbf{u}^{Q_n},
\label{convex}
\end{equation}
with
\begin{equation}
\lambda = \frac{1}{2} \left[1+\sin\left(\pi \left(\sigma_n-\frac{1}{2} \frac{\sigma_{\max} + \sigma_{\min}}{\sigma_{\max} -\sigma_{\min}}\right)\right)\right].
\end{equation}

A major concern for any shock capturing method is to maintain conservation, which ensures the correct shock speeds are maintained \cite{lax1960systems}.
\revCom{
\begin{remark}[Conservation] In its current incarnation the multi-element SIAC filter \textit{does not} conserve the solution quantities, e.g., density, momentum and energy for the Euler equations. The unfiltered standard DGSEM conserves the solution variables up to machine precision, see e.g. \cite{cockburn2000development}. However, the application of the local SIAC filter after each time step is no longer globally conservative because we re-distribute solution data, e.g. the mass, by the filtering process within each element. Whereas the conservation errors for the global filter are introduced by the necessary large interpolation order $N\gg m$, we can easily assure $N\leq m$ for the local element-wise filtering. However, we run into a different problem because our global approximation is no longer a polynomial, but solely built from piecewise polynomial data. Thus, again, we introduce conservation errors in our approximation, which are usually small. We examine the size of these conservation errors in Sec.~\ref{sec:convAndConsStudy} and show for the considered test cases that the conservation loss does not affect the solutions.
\end{remark}
}

\subsection{Two-dimensional filter}

Next, we extend the one-dimensional local SIAC filter to higher spatial dimensions. For the filtering process we apply the same local smoothing matrix $\hat{\mathbf{\Phi}}$ as in the one-dimensional case in each spatial direction to the unfiltered element solutions. Conveniently, this is possible due to the tensor product ansatz of the DGSEM and the definition of the Dirac delta kernel \eqref{delta}. Just as in the previous section we first derive a global multi-dimensional filter and then modify it for the local multi-element case.

We begin again from the filtering assumption of \eqref{start}, and find for the piecewise polynomial solution $u$ that
\begin{equation}
\begin{aligned}
\tilde{u}(x,y,t) & = \int\limits_{x-\varepsilon}^{x+\varepsilon} \int\limits_{y-\varepsilon}^{y+\varepsilon} u(\tau,\varsigma,t) \delta^{m,k}_{\varepsilon} (x-\tau,y-\varsigma) \dtaudsigma\\
           & = \sum_{\ell=1}^{N_Q} \int\limits_{x-\varepsilon}^{x+\varepsilon} \int\limits_{y-\varepsilon}^{y+\varepsilon} u^{Q_\ell}(\tau,\varsigma,t) \delta (x-\tau,y-\varsigma) \dtaudsigma.
\end{aligned}
\end{equation}
where we define the multi-variable delta function to have the form 
\begin{equation}
\delta(x,y) := \delta^{m,k}_{\varepsilon} (x,y) = \delta^{m,k}_{\varepsilon} (x)\delta^{m,k}_{\varepsilon} (y) =: \delta(x)\delta(y).
\end{equation}
We focus on one LGL node, transform into the reference space and use the tensor product property to split the integrand and obtain
\begin{equation}
\begin{aligned}
\tilde{u}_{ij}^{Q_n}(t) & = \sum_{\ell=1}^{N_Q} \int\limits_{x_i-\varepsilon}^{x_i+\varepsilon} \int\limits_{y_j-\varepsilon}^{y_j+\varepsilon} u^{Q_\ell}(\tau,\varsigma,t) \delta(x_i-\tau,y_j-\varsigma) \dtaudsigma \\
			& = \sum_{\ell=1}^{N_Q} \varepsilon^2 \int\limits_{-1}^1 \int\limits_{-1}^1 u^{Q_\ell}(\varepsilon x+\xi_i,\varepsilon y + \eta_j,t) \delta(\varepsilon x, \varepsilon y) \dxdy \\
			& = \sum_{\ell=1}^{N_Q} \varepsilon^2 \int\limits_{-1}^1 \int\limits_{-1}^1 \sum_{k,l=0}^N \bar{\psi}_k(\varepsilon x+\xi_i) \bar{\psi}_l(\varepsilon y + \eta_j) u_{kl}^{Q_{\bar{\ell}}}(t) \delta(\varepsilon x) \delta(\varepsilon y) \dxdy \label{TPsplit} \\
				& = \sum_{\ell=1}^{N_Q} \sum_{k,l=0}^N \underbrace{\varepsilon \int\limits_{-1}^1 \bar{\psi}_k(\varepsilon x+\xi_i) \delta(\varepsilon x) \dx}_{=\Phi_{ik}} ~ \underbrace{ \varepsilon \int\limits_{-1}^1 \bar{\psi}_l(\varepsilon y + \eta_j) \delta(\varepsilon y) \,\mathrm{d}y}_{=\Phi_{jl}}  ~ u_{kl}^{Q_{\bar{\ell}}}(t).
\end{aligned}
\end{equation}
Here, the $\bar{\ell}$ points to the correct solution entry, which includes neighboring elements and is dependent on the storing data structure. The shifting of the evaluation points for the lagrange basis function is also hidden in the bar notation, i.e. 
\begin{equation}
\bar{\psi}(x) := \begin{cases} ~~\psi(x), & x\in\left[-1,1\right] \\ \psi(x-2), & x>1 \\ \psi(x+2), & x<-1 \end{cases} 
\end{equation}
From this definition of the filtering matrices it is possible to write the filtering matrix in a compact notation
\begin{equation}
   \tilde{\mathbf{u}}^n = \mathbf{\Phi} \mathbf{u}_{env}^n \mathbf{\Phi}^T,
	\label{2D-MV}
\end{equation}
with	
\begin{equation}
	 \mathbf{\Phi} = \begin{pmatrix} \mathbf{\Phi}^{n-1} & \mathbf{\Phi}^n & \mathbf{\Phi}^{n+1} \end{pmatrix} ~~~ \text{  and } ~~~ \mathbf{u}_{env}^n = \begin{pmatrix} \mathbf{u}^{n+N_{Q,x}-1} & \mathbf{u}^{n+N_{Q,x}} & \mathbf{u}^{n+N_{Q,x}+1} \\ \mathbf{u}^{n-1} & \mathbf{u}^n & \mathbf{u}^{n+1} \\ \mathbf{u}^{n-N_{Q,x}-1} & \mathbf{u}^{n-N_{Q,x}} & \mathbf{u}^{n-N_{Q,x}+1} \end{pmatrix},
\end{equation}
provided the elements are labeled from bottom-left to top-right and $N_{Q,x}$ denotes the number of elements in the $x$-direction. In this case, we design the smoothing matrix $\mathbf{\Phi} = \begin{pmatrix} \mathbf{\Phi}^{n-1} & \mathbf{\Phi}^n & \mathbf{\Phi}^{n+1} \end{pmatrix}$ exactly as in one spatial dimension.

We want the resulting shock capturing DG scheme to be as local as possible and implement the 2D multi-element SIAC filter in a way to reflect this goal. First, we define
\begin{equation}
\begin{pmatrix} \hat{\mathbf{u}}^{n+N_{Q,x}} \\ \hat{\mathbf{u}}^n \\ \hat{\mathbf{u}}^{n-N_{Q,x}} \end{pmatrix} := \mathbf{u}_{env}^n \mathbf{\Phi}^T = \begin{pmatrix} \Phi^{n-1} \mathbf{u}^{n+N_{Q,x}-1} +  \Phi^n \mathbf{u}^{n+N_{Q,x}} + \Phi^{n+1} \mathbf{u}^{n+N_{Q,x}-1} \\  \Phi^{n-1} \mathbf{u}^{n-1} + \Phi^n \mathbf{u}^{n} + \Phi^{n+1} \mathbf{u}^{n+1}  \\  \Phi^{n-1} \mathbf{u}^{n-N_{Q,x}-1} +  \Phi^n \mathbf{u}^{n-N_{Q,x}} + \Phi^{n+1} \mathbf{u}^{n-N_{Q,x}-1} \end{pmatrix},
\end{equation}
which is nothing more than the solution vector of the three considered adjacent cells filtered in the $x$-direction. To filter in the $y$-direction, we just apply the smoothing matrix $\mathbf{\Phi}$ from the left hand side and obtain an overall filtered solution $\tilde{\mathbf{u}}^n$ from \eqref{2D-MV}. The main advantage of this procedure is that the filtering procedure is done dimension by dimension. So for all elements $n=1,\ldots,N_Q$ we first filter in the $x$-direction to find $\hat{\mathbf{u}}^{n}$ with coupling only from the right and left neighbor cells. Next, we filter in the $y$-direction and compute the fully filtered solution $\tilde{\mathbf{u}}^n$ from the information stored in the intermediate array $\hat{\mathbf{u}}^{n}$ with coupling from the upper and lower neighbor elements. That is, we simply apply the one-dimensional filter twice for each grid point and, again, only need information from the direct neighbors. Additionally, we note, that this filtering procedure has no preferred direction such that the order of $x,y$ directions makes no difference.

\section{Numerical tests}\label{Sec:NumTest}

We verify the performance of the novel multi-element SIAC filter applied to a variety of two-dimensional shock tests for the Euler as well as ideal MHD equations. For all simulations, we consider two-dimensional Cartesian meshes discretized by uniform quadrilateral elements of equal size $\Delta x = \Delta y$. Further, we use the explicit five stage, fourth order Runge-Kutta method of Carpenter and Kennedy \cite{Kennedy1994} with an stable explicit time step to advance the DG approximation in time. The (adaptive) filtering procedure \eqref{2D-MV} is performed for all element solutions after each time step. We begin with the numerical validation for the two-dimensional Euler equations in Sec.~\ref{Sec:Euler}, where we also investigate the accuracy and conservation issues of the filter, before we apply it to more challenging shock tests for the ideal MHD equations in Sec.~\ref{Sec:MHD}.

\subsection{Euler tests}\label{Sec:Euler}
The two-dimensional Euler equations are described by a system of conservation laws, i.e.
\begin{equation}
U_t + F_x + G_y  = 0,
\label{Euler2D}
\end{equation}
with
\begin{equation}
U = \begin{pmatrix} \rho \\ \rho \vec{v} \\ \rho e \end{pmatrix}, ~~~ F = \begin{pmatrix} \rho v_1 \\ \rho v_1^2 + p \\ \rho v_1 v_2 \\ v_1 \left(\rho e + p\right) \end{pmatrix}, ~~~ G = \begin{pmatrix} \rho v_2 \\ \rho v_1 v_2 \\ \rho v_2^2 + p \\ v_2 \left(\rho e + p\right) \end{pmatrix}.
\label{EulerFluxes}
\end{equation}
Here, $\rho$, $\vec{v}=(v_1, v_2)^T$, and $e$ denote the density, two-dimensional velocity field and inner specific energy, respectively. We close the system by the perfect gas equation, which relates the inner energy and pressure, i.e.
\begin{equation}
p = (\gamma-1)\left(\rho e - \frac{1}{2}\rho\|\vec{v}\|^2\right),
\end{equation}
where $\gamma$ denotes the adiabatic coefficient.

In the following, we apply the DGSEM with the multi-element SIAC filter derived herein to the two-dimensional Euler equations \eqref{Euler2D} to first show the high-order convergence and investigate the conservations properties of the scheme. Thereupon, we consider several benchmark problems for the two-dimensional Euler equations in order to verify the shock capturing capabilities of the novel filtering strategy.

\subsubsection{Convergence and conservation studies}\label{sec:convAndConsStudy}
A substantial property of DG schemes is the high-order accuracy of the approximation for smooth solutions. Thus, we first consider an academic test case with a known analytical solution for the two-dimensional Euler equations, that allows us to compute numerical errors measured in a discrete $L_\infty$-norm. With the help of these error values for different discretization levels, we are able to compute the experimental order of convergence (EOC), which for the DGSEM is expected to agree with the theoretical order of $N\! +\! 1$ as the mesh is refined. 

The problem we use for convergence tests is defined in the domain $\Omega = [-1,1]^2$ with the initial conditions
\begin{equation}
\rho(x,y,0) = 1 + 0.3 \sin(2 \pi (x+y)), \quad v_1(x,y,0) = v_2(x,y,0) = p(x,y,0) = 1.
\end{equation}
We use periodic boundary conditions and set $\gamma=5/3$. The analytical solution of the 2D Euler equations using these initial conditions is
\begin{equation}
U(x,y,t) = \left( \begin{array}{c} 
1 + 0.3 \sin (2 \pi (x+y-2t)) \\ 
1 + 0.3 \sin (2 \pi (x+y-2t)) \\
1 + 0.3 \sin (2 \pi (x+y-2t)) \\  
0.5+0.15 \sin(2 \pi (x+y-2t))+ \frac{1}{\gamma -1} \end{array} \right).
\end{equation}
We run the simulation until the final time $T = 0.4$ and intentionally choose a high polynomial degree of $N=7$. Consequently, we select a small explicit time step where $\texttt{CFL} = 0.1$ to exclude errors from the time integration method. Further, we compute the errors of the approximation for different choices of $N_Q$ and calculate  the EOC by the maximum error $\epsilon_\infty$ of the density. For the DGSEM approximation without filtering we obtain the convergence results illustrated in Table~\ref{tab:table6}.
\begin{table}[!ht]
  \caption{Euler convergence test for $N=7, \texttt{CFL} = 0.1$ and $T = 0.4$ without filtering.}
  \centering
  \begin{tabular}{@{}llll@{}}
    \toprule
    $N_Q$ & \multicolumn{1}{c}{$\epsilon_\infty$} & $EOC_\infty$ &  \multicolumn{1}{c}{$\epsilon_\mathrm{cons}$} \\
    \midrule
   $1$ & $5.72 \cdot 10^{-3}$  & ---  & $4.4 \cdot 10^{-16}$ \\
   $2^2$ & $4.56 \cdot 10^{-5}$  & $\bm{6.97}$  & $1.3 \cdot  10^{-15}$  \\
   $4^2$ & $1.74 \cdot 10^{-7}$  & $\bm{8.03}$  & $3.1 \cdot 10^{-15}$ \\
   $8^2$ & $4.82 \cdot 10^{-10}$  & $\bm{8.50}$  & $8.9 \cdot 10^{-16}$ \\
   $16^2$ & $1.79 \cdot 10^{-12}$  & $\bm{8.07}$  &  $2.0 \cdot 10^{-14}$\\
   \bottomrule
  \end{tabular}
  \label{tab:table6}
\end{table}

The results in Table~\ref{tab:table6} confirm the expected theoretical order of convergence $N+1$ for the DGSEM approximation without filtering. Moreover, we state the overall conservation error of the density computed by
\begin{equation}
\epsilon_\mathrm{cons} = \left|\int_\Omega \rho_{\text{ex}}(x,y,T) \dxdy - \int_\Omega \rho_{\text{app}}(x,y,T)\dxdy\right|,
\label{conserror}
\end{equation}
which regard to machine precision and, thus, agree with the desired conservative nature of the DGSEM \cite{cockburn2000development}.

As we are interested in the errors and convergence rates of the filtered solution as well, we turn off the adaptive filtering and investigate the convergence rates and conservation errors of the filtered solution as demonstrated in Table~\ref{tab:table7} for the SIAC filter using a Dirac-delta approximation with one vanishing moment.  
\begin{table}[!ht]
  \caption{Euler convergence test for $N = 7, \texttt{CFL} = 0.1$ and $T = 0.4$ filtered with $(m,k) = (1,6), N_d = 0.8$.}
  \centering
  \begin{tabular}{@{}llll@{}}
    \toprule
    $N_Q$ & \multicolumn{1}{c}{$\epsilon_\infty$} & $EOC_\infty$ & \multicolumn{1}{c}{$\epsilon_\mathrm{cons}$}  \\
    \midrule
   $10^2$ & $7.39 \cdot 10^{-2}$  & ---  & $5.7 \cdot  10^{-6}$  \\
   $20^2$ & $3.97\cdot 10^{-2}$  & $\bm{0.90}$  & $4.6 \cdot 10^{-6}$ \\
   $40^2$ & $2.06 \cdot 10^{-2}$  & $\bm{0.95}$  & $1.8 \cdot 10^{-6}$ \\
   $80^2$ & $1.05 \cdot 10^{-2}$ &$ \bm{0.97}$  & $5.4 \cdot 10^{-7}$\\
   \bottomrule
  \end{tabular}
  \label{tab:table7}
\end{table}

Table \ref{tab:table7} reveals that the order of convergence drops to one. Further, we lose conservation as pointed out in the previous section. We now repeat the convergence test with $m = 3$ and $m = 5$, see Tables \ref{tab:table8} and \ref{tab:table9}. 
\begin{table}[!ht]
  \caption{Euler convergence test for $N = 7, \texttt{CFL} = 0.1$ and $T = 0.4$ filtered with $(m,k) = (3,6), N_d = 2.5$.}
  \centering
  \begin{tabular}{@{}llll@{}}
    \toprule
    $N_Q$ & \multicolumn{1}{c}{$\epsilon_\infty$} & $EOC_\infty$ & \multicolumn{1}{c}{$\epsilon_{\mathrm{cons}}$}  \\
    \midrule
   $2^2$ & $2.52 \cdot 10^{-2}$  & ---  & $1.1 \cdot 10^{-3}$ \\
   $4^2$ & $3.56 \cdot 10^{-3}$  & $\bm{2.82}$  & $3.8 \cdot  10^{-5}$  \\
   $8^2$ & $4.46\cdot 10^{-4}$  & $\bm{3.00}$  & $1.3 \cdot 10^{-6}$ \\
   $16^2$ & $5.56 \cdot 10^{-5}$  & $\bm{3.00}$  & $3.6 \cdot 10^{-8}$\\
   \bottomrule
  \end{tabular}
  \label{tab:table8}
\end{table}

\begin{table}[!ht]
  \caption{Euler convergence test for $N = 7, \texttt{CFL} = 0.1$ and $T = 0.4$ filtered with $(m,k) = (5,7), N_d = 4.5$.}
  \centering
  \begin{tabular}{@{}llll@{}}
    \toprule
    $N_Q$ & \multicolumn{1}{c}{$\epsilon_\infty$} & $EOC_\infty$ & \multicolumn{1}{c}{$\epsilon_{\mathrm{cons}}$}  \\
    \midrule
   $1$ & $1.14 \cdot 10^{-3}$  & ---  & $9.5 \cdot 10^{-5}$ \\
   $2^2$ & $4.35 \cdot 10^{-5}$  & $\bm{4.71}$  & $1.4 \cdot  10^{-6}$  \\
   $4^2$ & $1.35 \cdot 10^{-6}$  & $\bm{5.01}$  & $1.2 \cdot 10^{-8}$ \\
   $8^2$ & $4.21 \cdot 10^{-8}$  & $\bm{5.00}$  & $8.9 \cdot 10^{-11}$ \\
   \bottomrule
  \end{tabular}
  \label{tab:table9}
\end{table}
As the convergence tests show, the experimental order of convergence depends on the number of vanishing moments in the Dirac-delta approximation, i.e. $EOC \approx {\min}\left(m,N+1\right)$. While one vanishing moment leads to a smearing effect and significantly drops the accuracy in smooth areas, increasing the number of vanishing moments improves the accuracy again. We observe that the conservation error decreases with an increasing number of vanishing moments as well as with finer grid resolutions.

In conclusion, we expect the Dirac-delta filter to effectively distinguish between shocks and smooth areas when using enough vanishing moments. Because of the large drop of accuracy for less vanishing moments, smooth areas will be mistaken for shocks more often which results in significantly less accuracy in those areas. We alleviate  these issues with the adaptive filtering technique, so that small conservation errors are exclusively introduced in shocked regions. 

\subsubsection{Explosion problem}
We now consider the first two-dimensional Euler test case that inherits shocks. The explosion problem is defined on the domain $\Omega = [-1,1]^2$, whereas its initial conditions consist of a region inside a circle with the radius $r = 0.4$ centered at the origin and a region outside that circle, see e.g. \cite{wissink2018shock}. The primitive variable initial conditions are then defined by two states, i.e. 
\begin{equation}
\rho_{\text{in}} (x,y,0) = 1.0, \quad \rho_{\text{out}} (x,y,0) = 0.125, \quad\quad \quad  p_{\text{in}} (x,y,0)= 1.0, \quad  p_{\text{out}} (x,y,0)= 0.1. 
\end{equation}
Here, the inner state applies for every $(x,y)\in\Omega$ such that $(x^2+y^2) \leq r^2$ and the outside state applies otherwise. Further, the velocity vector is set to zero in the entire domain and we set $\gamma = 5/3$.

\begin{figure}[!ht]
\centering
\includegraphics[height=5.3cm]{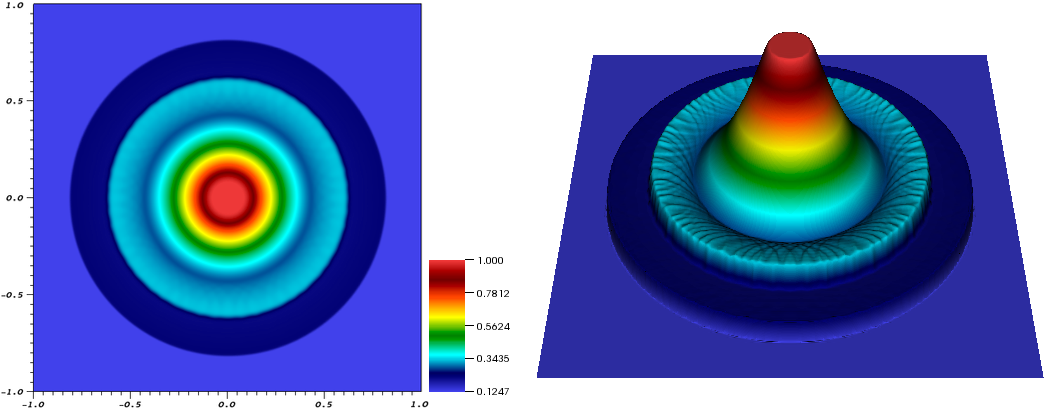}
\caption{Density of the 2D explosion problem at $T = 0.25$ for $N=7, N_Q = 80^2, \texttt{CFL} = 0.1$ filtered adaptively with $(m,k)=(1,6)$, $N_d=0.6, \sigma_{\min} = -7$ and $\sigma_{\max} = -3$.}
\label{Explosion_m1}
\end{figure}

\begin{figure}[!ht]
\centering
\includegraphics[height=5.3cm]{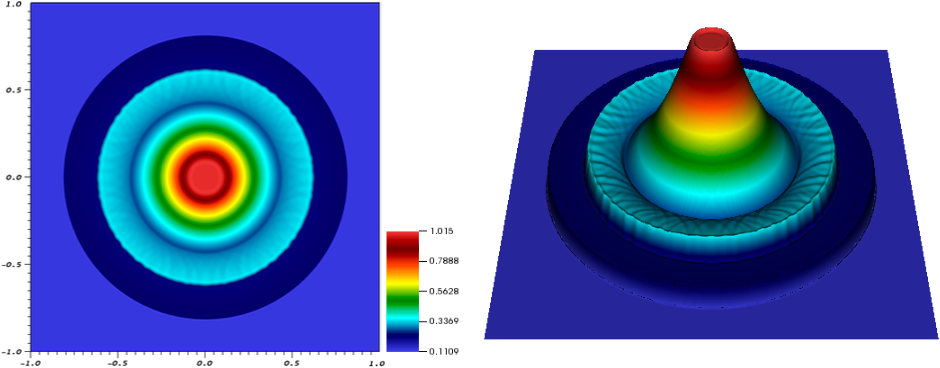}
\caption{Density of the 2D explosion problem at $T = 0.25$ for $N=7, N_Q = 80^2, \texttt{CFL} = 0.1$ filtered adaptively with $(m,k)=(3,6)$, $N_d=2.5, \sigma_{\min} = -8$ and $\sigma_{\max} = -5$.}
\label{Explosion_m3}
\end{figure}

The following plots in Figures \ref{Explosion_m1} and \ref{Explosion_m3} show the approximation of the density for this problem at the final time $T=0.25$ with $\texttt{CFL} = 0.1$ and a polynomial degree of $N=7$ on $N_Q = 80^2$ elements.

For the simulation result illustrated in Fig.~\ref{Explosion_m1}, one vanishing moment ($m=1$) and six continuous derivatives at the endpoints ($k=6$) of the Dirac-delta approximation are chosen. The support width in \eqref{support} is calculated with $N_d = 0.6$. Furthermore, the adaptive filtering with convex blending \eqref{convex} is used with a minimum tolerance of $10^{-7}$ and a maximum tolerance of $10^{-3}$. Comparing the results to other configurations of the filter shows that the narrow domain width and the chosen tolerances effectively limit the smearing effect despite the small number of vanishing moments. The conservation error is $\epsilon_\mathrm{cons} = 6.4 \cdot 10^{-4}$.

In contrast, Fig.~\ref{Explosion_m3} uses three vanishing moments ($m=3$), which requires a larger support width defined by $N_d = 2.5$ in \eqref{support}. The corresponding conservation error is $\epsilon_\mathrm{cons} = 2.0 \cdot 10^{-5}$.

\begin{figure}[!ht]
\centering
\includegraphics[scale=0.75]{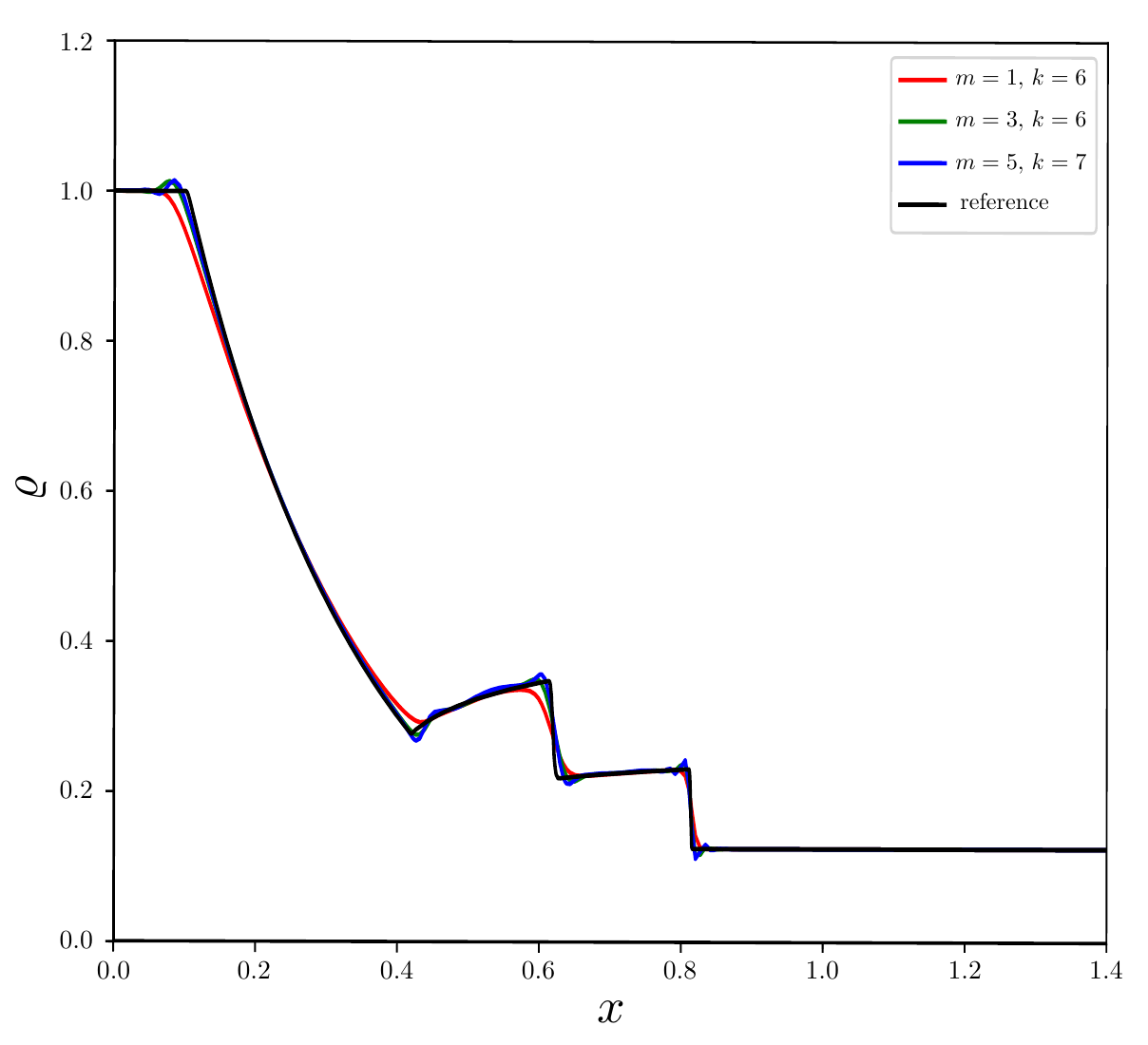}
\caption{Density slices of the 2D explosion problem at $x = y$, $T = 0.25$ for $N = 7$, $N_Q = 80^2$, \texttt{CFL} = 0.1.}
\label{Explosion_ref}
\end{figure}

Both results reveal that the shocks are effectively regularized, even though small overshoots are observable. For this problem, we find that one vanishing moment is enough to obtain acceptable results if the tolerance of the adaptive filtering operation is adjusted. However, the three moment adaptive SIAC filter produces more accurate shock profiles, as can be seen in the slicing picture Fig.~\ref{Explosion_ref} that compares three filtered solutions against a reference solution \cite{wissink2018shock}. The reference solution is computed by the publicly available high performance application code \texttt{FLASH} (http://flash.uchicago.edu/site/flashcode/) with a second order MUSCL-Hancock finite volume method (see e.g. \cite{Waagan2009}) on $2048 \times 2048$ elements. We observe that the filter with one vanishing moment is quite dissipative, whereas the filtered approximations with higher moments produce small overshoots at the shock and rarefraction. Nonetheless, all three configurations are stable, nearly oscillation-free and close to the analytical profile.

\subsubsection{Four state Riemann test}
The next Euler shock test is described by the 2D-Riemann problem, see \cite{kurganov2002solution} and \cite{lax1998solution}. We consider the domain $\Omega = [0,1]^2$ with outflow boundary conditions. The initial conditions are defined by four different states in four quadrants:
\begin{description}
\item[\textbf{Top left:}] $\Omega_{t\ell} = [0,0.5)\times(0.5,1]$.\\[-0.35cm]
\item[\textbf{Bottom left:}] $\Omega_{b\ell} = [0,0.5)\times(0,0.5]$.\\[-0.35cm]
\item[\textbf{Top right:}] $\Omega_{tr} = (0.5,1]\times(0.5,1]$.\\[-0.35cm]
\item[\textbf{Bottom right:}]  $\Omega_{b r} = (0.5,1]\times[0,0.5)$.
\end{description}
The four primitive variable states of the initial conditions are then assigned to these four quadrants:
\begin{align*}
(\rho, v_1, v_2, p)(x,y,0) = 
\begin{cases} 
(\rho_{t\ell}, v_{1,t\ell}, v_{2,t\ell}, p_{t\ell}), \qquad & \text{if } (x,y) \in \Omega_{t\ell}, \\
(\rho_{b\ell}, v_{1,b\ell}, v_{2,b\ell}, p_{b\ell}), \qquad & \text{if } (x,y) \in \Omega_{b\ell}, \\
(\rho_{tr}, v_{1,tr}, v_{2,tr}, p_{tr}), \qquad & \text{if } (x,y) \in \Omega_{tr}, \\
(\rho_{br}, v_{1,br}, v_{2,br}, p_{br}), \qquad & \text{if } (x,y) \in \Omega_{br}.
\end{cases}
\end{align*}
Depending on the particular choice of initial conditions, the simulation has to overcome at least one rarefaction wave, shock wave or contact wave. In \cite{kurganov2002solution}, 19 particular configurations for initial values are given. Because configurations 17 and 19 deal with all of the mentioned waves, we apply the Dirac-delta filter to them. The results of the following configurations are obtained by using a polynomial degree of $N = 7$ on a grid of $60^2$ elements with $\texttt{CFL} = 0.1$.

\paragraph{Configuration 17:}
The four initial conditions of configuration 17 are defined for the primitive variables as
\begin{align*}
&\rho_{t\ell} = 2, 			&& v_{1,t\ell} = 0, 	&& v_{2,t\ell} = -0.3, 		&& p_{t\ell} = 1 		\\
&\rho_{b\ell} = 1.0625, 	&& v_{1,b\ell} = 0,	&& v_{2,b\ell} = 0.2145, 	&& p_{b\ell}=0.4 	\\
&\rho_{tr} = 1, 				&& v_{1,tr} = 0, 	&& v_{2,tr} = -0.4, 			&& p_{tr} = 1 			\\
&\rho_{br} = 0.5197, 		&& v_{1,br} = 0, 	&& v_{2,br} = -1.1259, 		&& p_{br} = 0.4.
\end{align*} 

\begin{figure}[!ht]
\centering
\includegraphics[height=7.25cm]{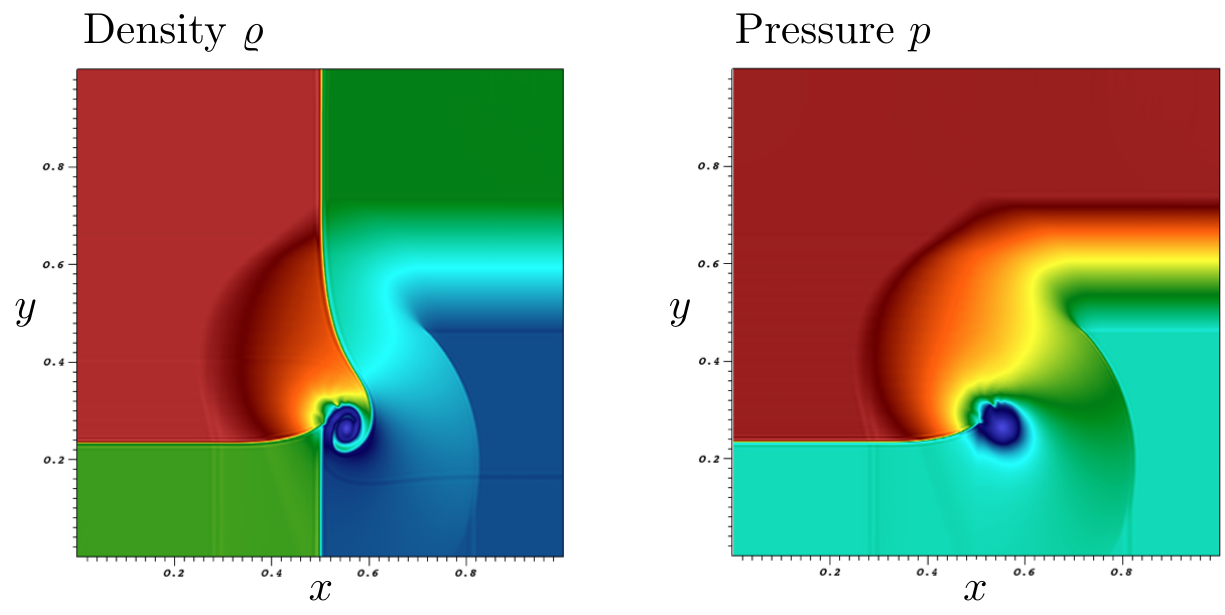}
\caption{Density (left) and pressure (right) of the four state Riemann test configuration 17 at $T=0.3$ for $N=7, N_Q = 60^2, \texttt{CFL}=0.1$ filtered adaptively with $(m,k)=(5,7), N_d=4.5, \sigma_\text{min} = -8$ and $\sigma_\text{max} = -3$.}
\label{conf17}
\end{figure}

We run the approximation to $T= 0.3$. The pseudo-color plot of the density in Fig.~\ref{conf17} is obtained by using the adaptive filter with a Dirac-delta approximation of $(m,k) = (5,7)$. The support width $\varepsilon$ is calculated by $N_d = 4.5$ in \eqref{support}, where the tolerances are set to $\sigma_\text{min} = -8$ and $\sigma_\text{max} = -3$. As illustrated in the plot, the filter ensures a successful regularization of shocks while keeping a high resolution of the vortex.

\paragraph{Configuration 19:}
The initial conditions for the primitive variables of configuration 19 are defined as
\begin{align*}
&\rho_{t\ell} = 2, 			&& v_{1,t\ell} = 0, 	&& v_{2,t\ell} = -0.3, 		&& p_{t\ell} = 1 		\\
&\rho_{b\ell} = 1.0625, 	&& v_{1,b\ell} = 0,	&& v_{2,b\ell} = 0.2145, 	&& p_{b\ell}=0.4 	\\
&\rho_{tr} = 1, 				&& v_{1,tr} = 0, 	&& v_{2,tr} = 0.3, 			&& p_{tr} = 1 			\\
&\rho_{br} = 0.5197, 		&& v_{1,br} = 0, 	&& v_{2,br} = -0.4259, 		&& p_{br} = 0.4.
\end{align*} 

\begin{figure}[!ht]
\centering
\includegraphics[height=7.25cm]{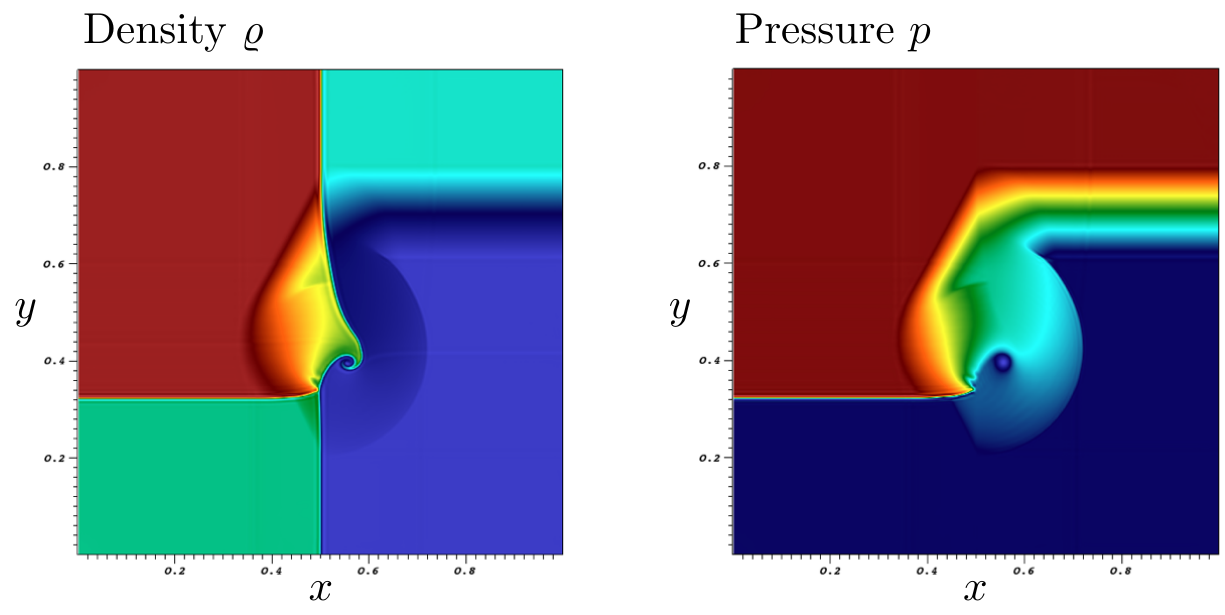}
\caption{DGSEM approximation of density $\rho$, left, and pressure $p$, right, adaptively filtered after every time step at $t = 0.3$ for $(m,k)=(5,7)$, $N_d=4.5, \sigma_\text{min} = 10^{-8}$ and $\sigma_\text{max} = 10^{-3}$.}
\label{conf19}
\end{figure}

Again, we filter adaptively after every time step with a Dirac-delta approximation defined by $m = 5, k=7$ and a support width calculated using $N_d = 4.5$. Figure \ref{conf19} demonstrates the approximation for the density and the pressure at the final time $T=0.3$, where the adaptive filtering again has minimum and maximum tolerances of $\sigma_\text{min} = -8$ and $\sigma_\text{max} = -3$. As before, the filter effectively regularizes shocked areas. Additionally, the adaptive filtering step ensures high resolution of the curly area in the center of the domain.

\subsubsection{Double Mach reflection}
As a final Euler test case, we consider the double Mach reflection problem, see e.g. \cite{shi2003resolution}. This challenging problem is defined on the domain $\Omega = [0,3.25]\times [0,1]$ and involves both, strong shock interactions and different boundary conditions. An initial shock separates a left and right state defined by the primitive variables, i.e. 
\begin{align*}
\rho_{L} (x,y,0)&= 8, & \rho_{R} (x,y,0)&= 1.4, \\
v_{1,L} (x,y,0)&= 8.25\cdot \frac{\pi}{6}, & v_{1,R} (x,y,0)&= 0, \\
v_{2,L} (x,y,0)&= -8.25\cdot \frac{\pi}{6}, & v_{2,R}(x,y,0) &= 0 ,\\
p_{L} (x,y,0)&= 116.5, & p_{R} (x,y,0)&= 1. 
\end{align*}
In particular, the shock is initially located at $x = 1/6, y=0$ along a linear slope, which includes an angle of $\alpha = \pi/3$ with the $x$-axis, defining the initial conditions
\begin{align*}
(\rho, v_1, v_2, p) = 
\begin{cases}
(\rho_L, v_{1,L}, v_{2,L}, p_L), \quad &\text{if } x <\frac{1}{6}+\frac{y}{\sqrt{3}}, \\
(\rho_R, v_{1,R}, v_{2,R}, p_R), \quad &\text{if } x \geq \frac{1}{6}+\frac{y}{\sqrt{3}}.
\end{cases}
\end{align*}

Thus, the boundary conditions at the bottom for $x < 1/6$ as well as at the left and upper domain boundaries are set to Dirichlet-values all through the simulation. Especially, the latter have to be implemented correctly by an analytical shock profile, since the shock moves forward and its position is determined by the following function
\begin{align*}
s(t) = \frac{1}{6} + \frac{1+20t}{\sqrt{3}}.
\end{align*}
Using this function, the values at the top boundary are set to
\begin{align*}
(\rho, v_1, v_2, p) = 
\begin{cases}
(\rho_L, v_{1,L}, v_{2,L}, p_L), \quad &\text{if } x <s(t), \\
(\rho_R, v_{1,R}, v_{2,R}, p_R), \quad &\text{if } x \geq s(t)
\end{cases}
\end{align*}
Furthermore, the right boundary conditions are outflow and the bottom is considered a reflecting wall for $x > 1/6$. Hence, in the filtering process the according ghost cell values are set constantly either to the analytical solution for Dirichlet or to the last interior value for outflow and reflection.

\begin{figure}[!ht]
\centering
\includegraphics[height=4.75cm,trim=12 10 10 10,clip]{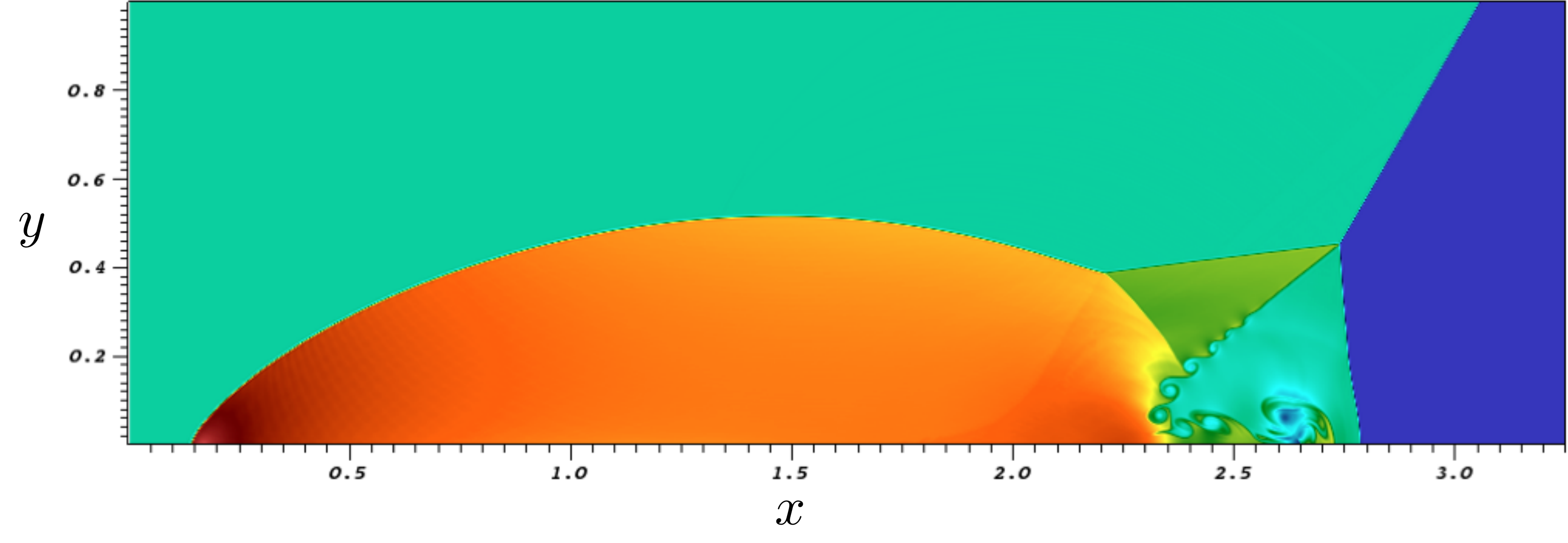}
\caption{Density of the double Mach reflection at $T = 0.2$ for $N=7,\texttt{CFL}=0.1$ on $325 \times 100$ elements filtered adaptively with $(m,k)=(3,6), N_d=2.5, \sigma_\text{min} = -7$ and $\sigma_\text{max} = -2$.}
\label{DoubleMach}
\end{figure}

In our simulation of this problem, we use a high grid resolution of $325 \times 100$ elements, a polynomial degree of $N=7$ and a CFL number of $\texttt{CFL} = 0.1$. We provide a density plot at the final time $T=0.2$ in Figure \ref{DoubleMach}, where we applied the adaptive multi-element SIAC filter with $m = 3, k = 6, N_d = 2.5$ and the according tolerances $\sigma_\text{min} = -6, \sigma_\text{max} = -2$. Compared to the previous shock tests, the number of vanishing moments had to be reduced in order to get smoother results. With these particular configurations, the approximation shows the curved reflected shock as well as the formation of a turbulent vortex, which compare well to the simulation results in the literature, e.g. \cite{sonntag2017}.

\subsection{Ideal MHD tests}\label{Sec:MHD}
Next, we consider the ideal magneto-hydrodynamic (MHD) equations, which can be expressed in a compact form as a system of conservation laws
\begin{equation}
U_t + F_x + G_y = 0
\label{MHD2D}
\end{equation}
with
\begin{equation}
U = \begin{pmatrix} \rho \\ \rho \vec{v} \\ \rho e \\ \vec{B} \end{pmatrix}, ~~~ F = \begin{pmatrix} \rho v_1 \\ \rho v_1^2 + p + \frac{1}{2}\|\vec{B}\|^2 - B_1^2 \\ \rho v_1v_2 - B_1B_2 \\ \rho v_1v_3-B_1B_3 \\ v_1\left[\rho e + p + \frac{1}{2} \|\vec{B}\|^2\right] - B_1(\vec{v}\cdot\vec{B}) \\ 0 \\ v_1B_2-v_2B_1 \\ v_1B_3-v_3B_1 \end{pmatrix} , ~~~ G = \begin{pmatrix} \rho v_2 \\ \rho v_1v_2 - B_1B_2 \\ \rho v_2^2 + p + \frac{1}{2}\|\vec{B}\|^2 - B_2^2 \\ \rho v_2v_3 - B_2B_3 \\ v_2\left[\rho e + p + \frac{1}{2} \|\vec{B}\|^2\right] - B_2(\vec{v}\cdot\vec{B}) \\ v_2B_1-v_1B_2 \\ 0 \\ v_2B_3-v_3B_2 \end{pmatrix} 
\label{MHDFluxes}
\end{equation}
Here, $\vec{v}=(v_1,v_2,v_3)^T$ denotes the velocity field and $\vec{B}=(B_1,B_2,B_3)^T$ the magnetic field components. As for the Euler equations, the system is closed by the perfect gas equation, which for the ideal MHD equations reads
\begin{equation}
p = (\gamma-1)\left(\rho e - \frac{1}{2}\rho\|\vec{v}\|^2 - \frac{1}{2}\|\vec{B}\|^2\right),
\end{equation}
where $\gamma$ again denotes the adiabatic coefficient. 

For the ideal MHD equations it is well known, that magnetic monopoles are not observed in nature, but due to round-off errors, numerically $\nabla \cdot \vec{B} = 0$ is not necessarily guaranteed. Hence, to maintain the divergence-free condition of the magnetic field $\nabla \cdot \vec{B} = 0$ we implement a hyperbolic divergence cleaning method based on generalized Lagrangian multiplier \cite{Dedner2002}. 

In the following we demonstrate the performance of the local SIAC Dirac-delta filter for two common benchmark problems of the two-dimensional MHD simulations including (strong) shocks. 
  
\subsubsection{Orszag-Tang vortex}
The first test case describes the evolution of a turbulent plasma cloud, the well-known Orszag-Tang vortex \cite{altmann2012}, which is initialized in the periodic domain $\Omega = \left[0,1\right]^2$ as
\begin{align*}
\rho = & ~ \gamma \, \frac{5}{12\pi},  ~~~~~~~~~~~~\, v_1 = -\sin(2\pi y), ~~~~~~~~~~~~\, v_2 = \sin(2\pi x), \\
p = & ~ \frac{5}{12\pi}, ~~~~~~~~~~~~~~ B_1 = -\frac{1}{\sqrt{4\pi}} \sin(2\pi y), ~~~~~\, B_2 = \frac{1}{\sqrt{4\pi}} \sin(4\pi x).
\end{align*}
The other variables are initially set to zero and $\gamma=5/3$. Furthermore, we use $\texttt{CFL} = 0.5$, polynomials of degree $N\!=\! 5$ and $40 \times 40$ elements for the following simulations of this test case.

\begin{figure}[!ht]
\hspace{-0.1cm}\includegraphics[width=14.75cm]{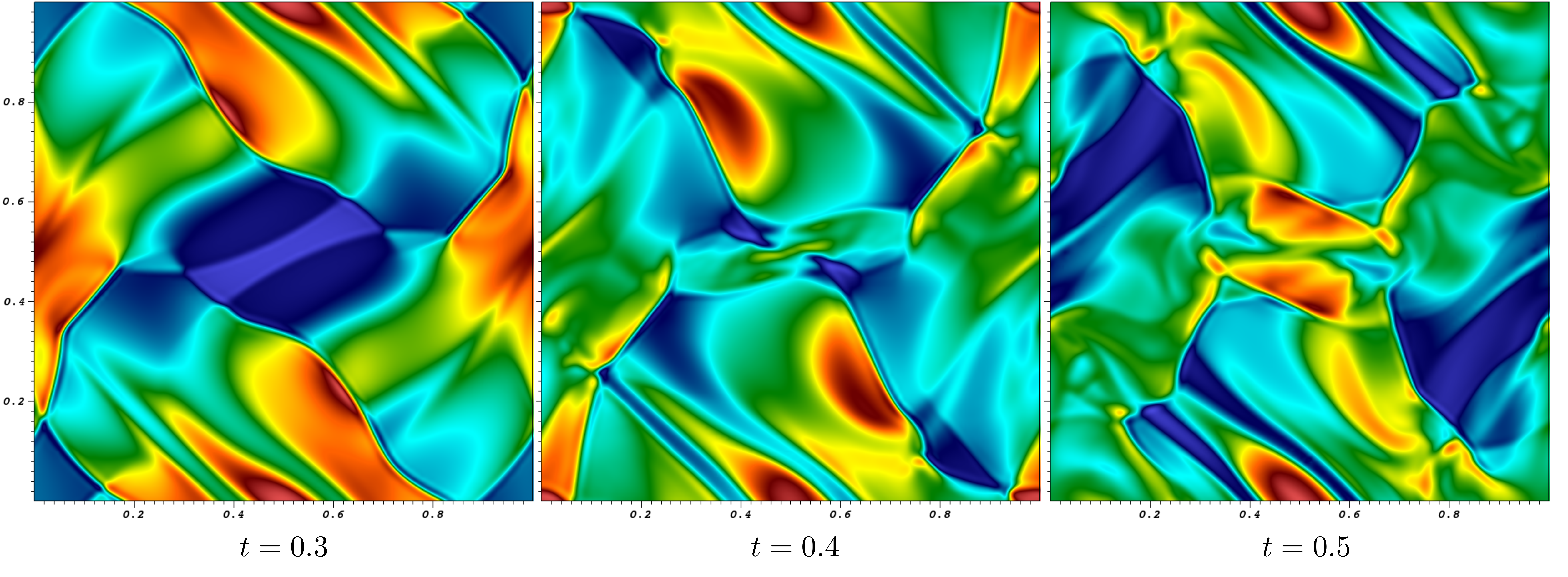}
\caption{Time evolution of the Orszag-Tang vortex density on $40 \times 40$ elements with $N\!=\! 5$ and $\texttt{CFL} = 0.5$ filtered adaptively with $m = 3, k = 8, \varepsilon = 1.4$ and $\sigma_{\min}=\sigma_{\max}=-8$.}
\label{Fig:OTVfiltered}
\end{figure}

\begin{figure}[!ht]
\hspace{-0.1cm}\includegraphics[width=14.75cm]{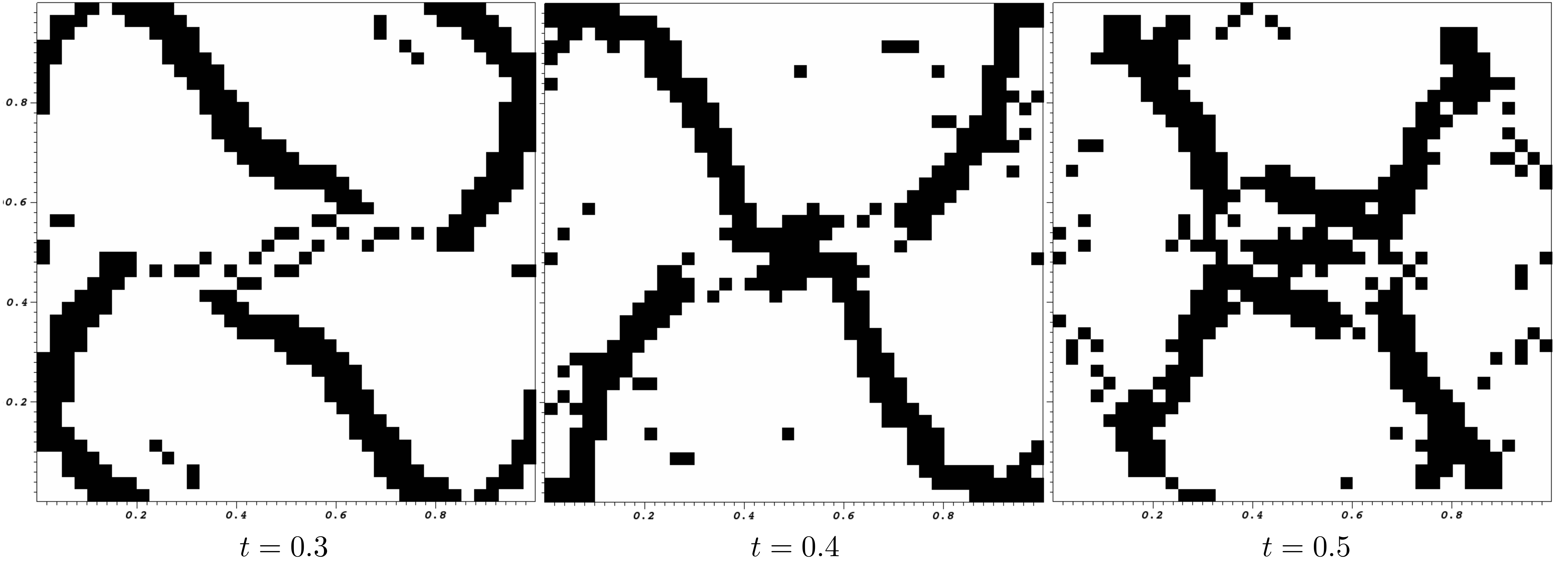}
\caption{Time evolution of the convex parameter $\lambda$ on $40 \times 40$ elements with $N\!=\! 5$ and $\texttt{CFL} = 0.5$ filtered adaptively with $m = 3, k = 8, \varepsilon = 1.4$ and $\sigma_{\min}=\sigma_{\max}=-8$.}
\label{Fig:OTV_lambda}
\end{figure}

We show the evolution of the density in the plots below (Fig.~\ref{Fig:OTVfiltered}) smoothed by the multi-element SIAC filter with $m = 3, k = 8$ and a fixed $\varepsilon = 1.4$. As for the Euler tests, we apply the filtering adaptively as shown in \eqref{convex} with the pressure as a shock indicator and $\sigma_{\min}=\sigma_{\max}=-8$. Further, we show the distribution of the cell-wise constant convex parameter $\lambda$ from \eqref{convex} at the same stages in Fig.~\ref{Fig:OTV_lambda}, which confirms the correct tracking of shocks as they evolve.

\begin{figure}[!ht]
\centering
\includegraphics[scale=0.585]{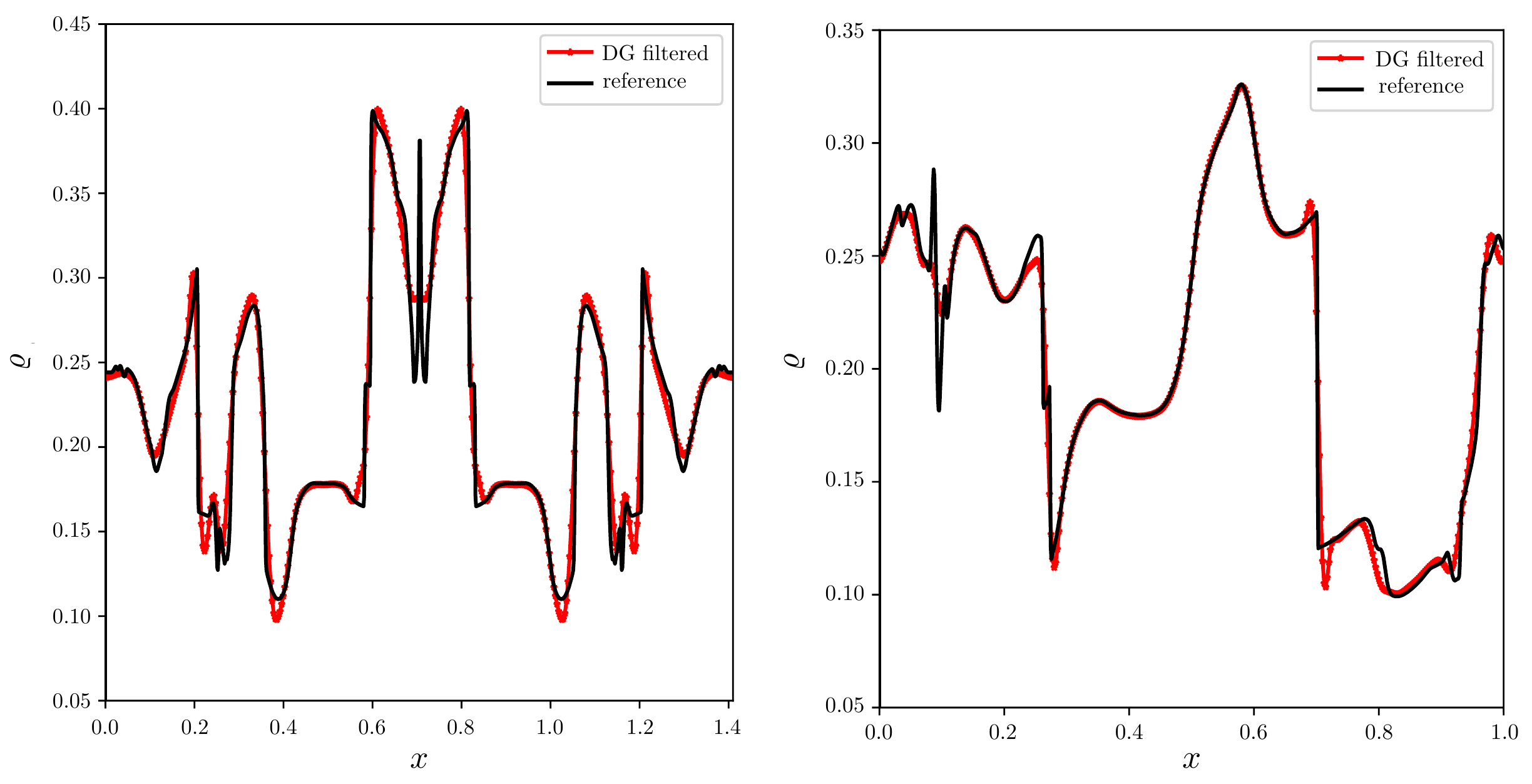}
\caption{Orszag-Tang-Vortex density slices at $x = y$ (left) and $y = 0.3$ (right) for $T = 0.5$, \texttt{CFL} = 0.5, $N = 5$ and $40 \times 40$ elements.}
\label{OTV_ref}
\end{figure}

In order to assess the performance of the two-dimensional Dirac-delta filter, we compare the simulation results to a reference solution obtained by the publicly available high performance application code \texttt{FLASH} (\url{http://flash.uchicago.edu/site/flashcode/}). Particularly, this highly resolved reference solution is computed by a second order MUSCL-Hancock finite volume method (see e.g. \cite{Waagan2009}) on $1024 \times 1024$ elements. Further, we apply the filter once more to the Orszag-Tang vortex with different parameters, i.e. we choose $m = 1, k = 6, \varepsilon = 1.6, \sigma_{\min}=-9$ and $\sigma_{\max}=-6$. In Fig.~\ref{OTV_ref} we provide two slices of the Orszag-Tang vortex at the final time $T=0.5$, in which we cut through the density distribution at $x=y$ (left) and $y=0.3$ (right) to compare the profiles obtained by the multi-element SIAC filter against a reference solution.

We see that the oscillations are smoothed out by the multi-element SIAC filter and the approximation matches the reference solution quite well, but the filtering technique still produces little overshoots at shocks. However, taking into account the coarse resolution of $40 \times 40$ elements, the filter generates reasonable approximations for this shock test, as it stabilizes the approximation and regularizes it against spurious oscillations.

\subsubsection{Magnetic rotor}
The second test case describes a rotating dense circle in a static fluid, that generates strong circular shock waves \cite{altmann2012}. In general this benchmark problem is defined in the same periodic domain $\Omega = \left[0,1\right]^2$ by the radius $r = \sqrt{(x-0.5)^2+(y-0.5)^2}$ and the slope $s=\frac{r_1-r}{r_1-r_0}$. The initial primitive variables for the magnetic rotor are stated in Table~\ref{tab:Rotor_states}, where the unlisted quantities are initially zero in the entire domain and $\gamma=1.4$.

\begin{table}[htbp!]
  \caption{Initial primitive states for the magnetic rotor test.}
  \centering
  \begin{tabular}{@{}lccccc@{}}
    \toprule
    Location & $\rho$ 	& $v_1$	& $v_2$	& $p$ 	& $B_1$	\\\midrule
   $r < r_0 $ & $10$	& $\dfrac{u_0}{r_0}\left(\frac{1}{2}-y\right)$	& $\dfrac{u_0}{r_0}\left(x-\frac{1}{2}\right)$	& $ ~ 1 ~$	&$\dfrac{5}{\sqrt{4\pi}}$\\[0.4cm]
   $r_0\leq r \leq r_1$  & $1+9s$	& $\dfrac{s u_0}{r_0}\left(\frac{1}{2}-y\right)$	& $\dfrac{s u_0}{r_0}\left(x-\frac{1}{2}\right)$& $1$	&$\dfrac{5}{\sqrt{4\pi}}$\\[0.4cm]
   $r>r_1$ & $1$	& $0$	& $0$& $1$	& $\dfrac{5}{\sqrt{4\pi}}$\\
   \bottomrule
   \end{tabular}
   \label{tab:Rotor_states}
\end{table}

\begin{figure}[!ht]
\includegraphics[width=14.7cm]{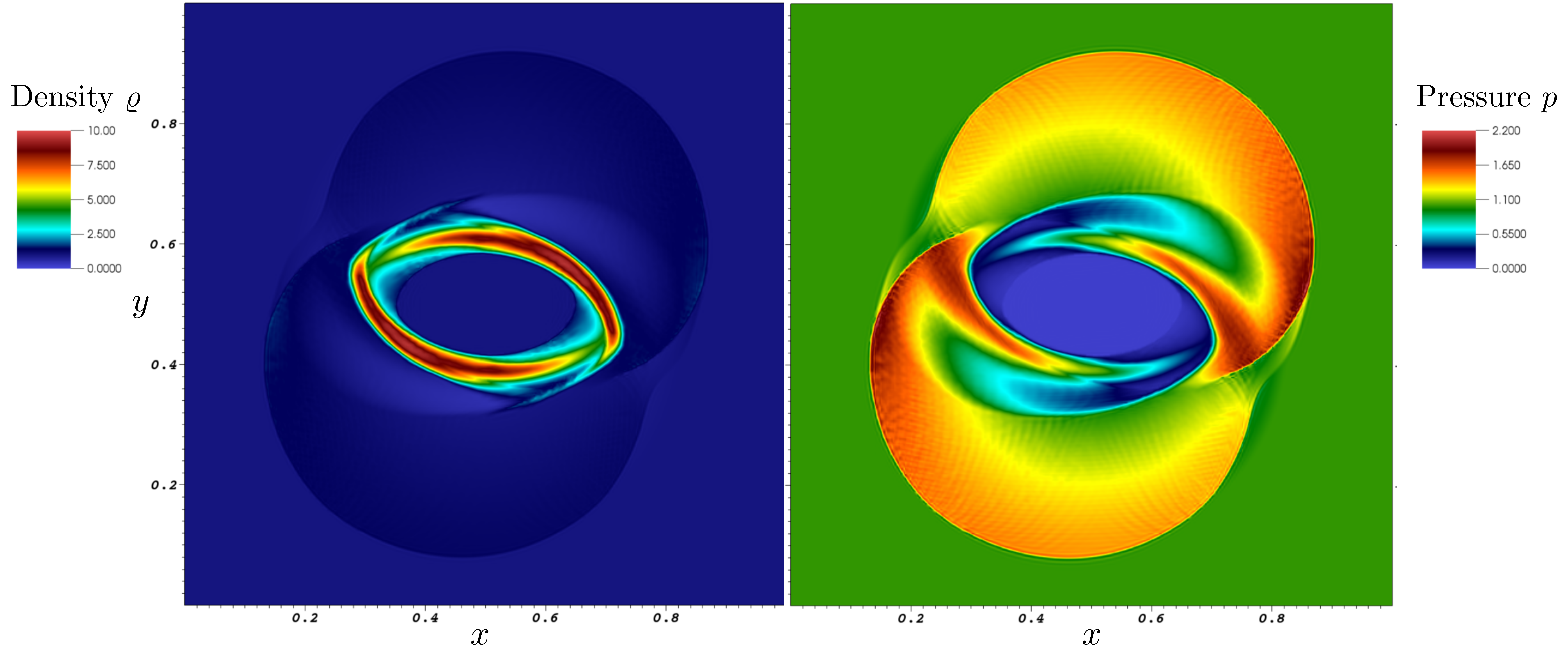}
\caption{Magnetic rotor density (left) and pressure (right) at $T=0.15$ for $\texttt{CFL}= 0.5, N\!=\! 4$ on $100 \times 100$ elements filtered adaptively with $(m,k)=(1,5), \varepsilon = 1.4, \sigma_\text{min} = -9$ and $\sigma_\text{max} = -6$.}
\label{rotor}
\end{figure}

In our simulations we define $r_0=0.1, r_1=0.115$ and $u_0=2$. We use $\texttt{CFL} = 0.5$, a polynomial degree of $N\!=\! 4$ and $100\times 100$ elements. We show the density and pressure at $T=0.15$ in Figure \ref{rotor}. Due to the strong circular shocks combined with the high-order DG approximation this test case is extremely sensitive and unstable. Therefore, we apply the SIAC filtering matrix constructed by a Dirac-delta kernel with only one vanishing moment $m = 1$ and $k = 5$. Again, we smooth the approximation adaptively with the density as a shock indicator, $\sigma_{\min}=-9, \sigma_{\max}=-6$ and a fixed $\varepsilon = 1.4$. 

In Figure \ref{rotor} we see, that the local SIAC filter performs well in terms of stabilizing the approximation and regularizing oscillatory regions, but is polluted by small mesh artifacts, which are particularly visible at the generated Alfvén waves in the pressure profile. Nonetheless, the novel filter produces reasonable approximations even for such a challenging test problem including strong circular shock waves.

\section{Conclusion}\label{Sec:Conc}
We have presented a novel shock capturing technique for DG methods based on multi-element SIAC filtering. In particular, we have first introduced the DGSEM on two-dimensional Cartesian meshes, before we have derived the local filter matrix from the global SIAC filtering approach constructed by approximations of Dirac-delta kernels. Moreover, we have designed an adaptive filtering strategy and extended the overall method to higher spatial dimensions. Finally, we have verified the applicability of the filter to a variety of challenging shock problems for both, the two-dimensional Euler and the ideal MHD equations.

We have demonstrated that the multi-element SIAC filter indeed performs well in terms of regularizing oscillatory regions and stabilizing the approximation. However, the main issue of the constructed filter is the introduction of spurious, albeit small, conservation errors, which might be avoided by projections onto discontinuous basis functions across cell interfaces. This, together with  the investigation of the entropic properties as well as the extension of the local SIAC filter to curvilinear elements are future research projects.    

\section*{Acknowledgements}
Gregor Gassner has been supported by the European Research Council (ERC) under the European Union's Eights Framework Program Horizon 2020 with the research project Extreme, ERC grant agreement no. 714487. Professor Jacobs acknowledges funding from NSF-DMS 1115705.

\bibliographystyle{plain}
\bibliography{references}

\end{document}